\documentclass[11pt]{article}
\usepackage[utf8]{inputenc}
\usepackage[margin=1in]{geometry}
\usepackage{amsmath,amssymb,amsthm,mathtools,mathrsfs}
\usepackage{enumitem}
\usepackage{hyperref}
\usepackage{booktabs}
\usepackage{array}
\usepackage{microtype}
\usepackage[T1]{fontenc}
\usepackage{lmodern}
\usepackage{tikz-cd}
\usepackage{titlesec}

\hypersetup{colorlinks=true,linkcolor=blue,citecolor=blue,urlcolor=blue}

\setlist{topsep=4pt,itemsep=1pt,parsep=0pt,partopsep=0pt}
\titlespacing*{\section}{0pt}{2.4ex plus.8ex minus.3ex}{.9ex}
\titlespacing*{\subsection}{0pt}{1.8ex plus.6ex minus.2ex}{.6ex}

\newcommand{\CC}{\mathbb C}
\newcommand{\PP}{\mathbb P}
\newcommand{\Aff}{\mathbb A}
\newcommand{\OO}{\mathcal O}
\newcommand{\Gr}{\operatorname{Gr}}
\newcommand{\Jac}{\operatorname{Jac}}
\newcommand{\Span}{\operatorname{Span}}
\newcommand{\rk}{\operatorname{rk}}
\newcommand{\im}{\operatorname{im}}
\newcommand{\coker}{\operatorname{coker}}
\newcommand{\tors}{\operatorname{tors}}
\newcommand{\Hom}{\operatorname{Hom}}
\newcommand{\Ext}{\operatorname{Ext}}

\newcommand{\Sym}{\operatorname{Sym}}
\newcommand{\ev}{\operatorname{ev}}
\newcommand{\drk}{\operatorname{drk}}
\newcommand{\divisor}{\operatorname{div}}

\newcommand{\GL}{\operatorname{GL}}

\newtheoremstyle{compactdefinition}
  {5pt plus 1pt minus 1pt}{5pt plus 1pt minus 1pt}
  {\normalfont}{}{\bfseries}{.}{.5em}{}
\newtheoremstyle{compactplain}
  {5pt plus 1pt minus 1pt}{5pt plus 1pt minus 1pt}
  {\itshape}{}{\bfseries}{.}{.5em}{}

\theoremstyle{compactdefinition}
\newtheorem{definition}{Definition}[section]
\newtheorem{example}[definition]{Example}

\newtheorem{remark}[definition]{Remark}
\theoremstyle{compactplain}
\newtheorem{theorem}[definition]{Theorem}
\newtheorem{proposition}[definition]{Proposition}
\newtheorem{lemma}[definition]{Lemma}
\newtheorem{corollary}[definition]{Corollary}
\newtheorem{fact}[definition]{Fact}
\newtheorem*{theoremA}{Theorem A}
\newtheorem*{theoremB}{Theorem B}
\newtheorem*{theoremC}{Theorem C}
\newtheorem*{theoremD}{Theorem D}
\newtheorem*{theoremE}{Theorem E}
\newtheorem*{theoremF}{Theorem F}

\title{Geometry of $m$-direct-sum morphisms}
\author{Duo Li, Yanjie Li and Francesco Russo}
\date{}

\begin{document}
\maketitle

\begin{abstract}
Let \(V\) be a complex vector space of dimension \(mk\) ($m\ge2$).
A morphism from a positive-dimensional projective variety
\(f\colon X\to\operatorname{Gr}(k,V)\) is \(m\)-direct-sum if
\(V=U_{x_1}\oplus\cdots\oplus U_{x_m}\) for every \(m\) pairwise distinct
points $x_1,\dots,x_m$ on $X$, where \(U_x\subset V\) is the subspace represented by \(f(x)\).  We show that $X$  is necessarily a smooth rational curve. When the so-called drill $d(f)$ of $f$ is $1$ or when the rank $k$ is at most $2$, we give  complete classifications of $f$ up to the automorphisms of $\operatorname{Gr}(k,V)$. As applications, we classify all so-called skew scrolls in $\PP^3$,  which completely answers a question of Day. Then we show that if \(X\subset\mathbb P(V)\) is a smooth,  projective, nondegenerate  subvariety whose higher-order Terracini locus is empty, then $X$ is the rational normal curve, which strengthens the critical Terracini rigidity discovered by Ballico and Chiantini. By using the same machinery, we also prove a dual version of Terracini rigidity. Finally, we show that $m$-unisolvent systems form an abelian category of finite length and that $7$ is the least diagonal exponent of a simple rank-three $2$-unisolvent system.
\end{abstract}

\tableofcontents

\section*{Introduction}
\phantomsection
\addcontentsline{toc}{section}{Introduction}
In this article, all the varieties are assumed to be integral varieties over the field of complex numbers.
Rational normal curves  can be characterized by the following property.
 If $f:X\to\PP(V)$ is a nonconstant morphism from a smooth projective variety, where
$\dim V\text{ is }m$ and the lines $\ell_x\subset V$ represented by $f(x)$ satisfy
$
 V=\ell_{x_1}\oplus\cdots\oplus\ell_{x_m}
$
for every $m$ pairwise distinct points $x_1,\ldots,x_m\in X$, then
$X\text{ is isomorphic to }\PP^1$.  Up to a linear change of coordinates, $f$ is the Veronese morphism
$\nu_{m-1}$.  Motivated by this characterization, we propose a higher-rank analogue.  Let $V$ be a complex vector space of
dimension $mk$ and  let $X$ be a  projective variety of positive dimension. \textbf{\textit{We do not assume a priori that $X$ is smooth}}. Assume
$
 f:X\longrightarrow\Gr(k,V),\quad x\longmapsto U_x\subset V,
$
is a morphism.  We say that $f$ has the \emph{$m$-direct-sum property} if
$
 V\text{ is decomposed as }U_{x_1}\oplus\cdots\oplus U_{x_m}
$
for every ordered $m$-tuple of pairwise distinct points of $X$.  So we may ask a natural question: For any fixed $k$,
can we classify all $m$-direct-sum morphisms up to the automorphisms of $\Gr(k,V)$?

The direct-sum condition is closely related to the theory of
$k$-regular maps.  Recall that a map to projective space is
projectively $m$-regular if the images of every $m$ distinct points
are projectively independent.  Buczy\'nski, Januszkiewicz, Jelisiejew
and Micha{\l}ek developed an algebro-geometric approach to the
construction of such maps using finite schemes, punctual Hilbert
schemes and areoles \cite{BJJM}; see also \cite{MM} for related
constructions and interpolation spaces.  More recently, Jelisiejew
and Keneshlou extended this framework to Grassmannians
\cite{JK}: a morphism
$  X\longrightarrow \operatorname{Gr}(k,V)$
is $m$-regular in their terminology if the $k$-planes corresponding
to any $m$ distinct points span a subspace of dimension $mk$.
Thus an $m$-direct-sum morphism in the present paper is precisely an
algebraic Grassmannian $m$-regular morphism in the critical ambient
dimension
$      \dim V=mk.$
The emphasis here is complementary to the construction problem
studied in \cite{BJJM,JK}.  Rather than seeking $m$-regular maps in
ambient spaces of small dimension, we fix the smallest possible
ambient dimension $mk$ and study the resulting rigidity and
classification problem for projective algebraic sources.  In this
critical situation projectivity is highly restrictive: we prove that
every positive-dimensional source is necessarily a smooth rational
curve, and we obtain complete classifications when the drill is one
or when the Grassmannian rank is at most two.  The Hilbert-scheme
compactification also acquires a quantitative form in this setting.
If $(E,W)$ is the coherent system associated with an $m$-direct-sum
curve and $e$ is its diagonal exponent, then the universal evaluation
over $\operatorname{Hilb}^m(\mathbf P^1)$ satisfies
$\operatorname{div}(\det\epsilon_{E,W})
      =\frac{e-k}{2}\Delta_m .$
Consequently, in the terminology of \cite{JK}, the morphism is
strongly $m$-regular if and only if $e=k$, equivalently if and only if
it is the split Vandermonde model.

The characterization of $m$-direct-sum curves can answer two concrete rigidity questions in projective
geometry.  The first is a problem of Day concerning skew scrolls in
$\PP^3$.  Let $C$ be a smooth connected projective curve and let
$\iota_1,\iota_2:C\hookrightarrow\PP^3$ be embeddings with
$\iota_1(x)\neq\iota_2(x)$ for every $x\in C$.  Joining corresponding points
gives a one-parameter family of lines
$
 \ell_x=\langle\iota_1(x),\iota_2(x)\rangle\subset\PP^3.
$
Its bidegree is
$\bigl(\deg\iota_1^*\OO_{\PP^3}(1),\deg\iota_2^*\OO_{\PP^3}(1)\bigr)$.
The resulting scroll is called \emph{skew} if we have
$\ell_x\cap\ell_y=\varnothing$ for all $x\neq y$.  Day proved that the
parameter curve is rational and, under a
genericity hypothesis, that the bidegree is $(1,1)$. He also constructed a skew scroll of bidegree $(2,2)$ and asked whether skew scrolls of higher
bidegrees exist (see
\cite[Lemma~3.11, Theorem~3.13, and Example~3.14]{Day2025}).

The second question is the critical Terracini problem.  Let
$X\subset\PP(V)$ be a smooth, connected and nondegenerate subvariety of dimension $n$. Assume
$\dim V\text{ is }m(n+1)$.  For $x\in X$, let $\widehat T_xX\subset V$ denote the
affine tangent space.  The emptiness of the reduced critical Terracini locus is equivalent to
$
 \widehat T_{x_1}X+\cdots+\widehat T_{x_m}X=V
$
for every $m$ pairwise distinct points.    Equivalently, the
tangent Gauss morphism
$
 \gamma_X:X\longrightarrow\Gr(n+1,V),\quad
 x\longmapsto\widehat T_xX,
$
is $m$-direct-sum.  Ballico and Chiantini proved that, for a curve in
$\PP^{2m-1}$, emptiness of the critical Terracini locus characterizes the
degree-$(2m-1)$ rational normal curve
\cite[Proposition~7]{BallicoChiantini2024}; see also
\cite{BallicoChiantini2021} for the general theory of Terracini loci.

The above two questions can both be reformulated in the language of $m$-direct-sum morphisms and we can attack them  by the machinery systematically developed in this article.

This article is organized as follows: In Section~\ref{sec:direct-sum-varieties},  we show that the source of any $m$-direct-sum morphism is a smooth rational curve in Theorem \ref{source} and prove an affine version of the characterization of rational normal curves. In Section~\ref{sec:perkinson}, we study
Perkinson's derived normal forms, especially the concept of the drill of a morphism into a Grassmannian variety.  We prove that if the drill of an $m$-direct-sum morphism is $1$, then the morphism is necessarily the osculating curve of a nondegenerate rational normal curve.
Section~\ref{sec:unisolvent} introduces the key invariant called  diagonal exponent.   Concretely, we associate an $m$-direct-sum curve
$f:\PP^1\to\Gr(k,V)$ with a special coherent system (see \cite{KingNewstead95} for the theory of coherent systems) on $\PP^1$ in the following way. Let
$S\text{ be }f^*\mathcal U$, where $\mathcal U$ is the tautological subbundle.  After
dualizing the inclusion $S\subset V\otimes\OO_{\PP^1}$, we obtain a generated
coherent system
$(E,W)$
 where  $E\text{ is }S^\vee\text{ with }W=V^\vee\subset H^0(\PP^1,E).
$
For every reduced divisor $D=x_1+\cdots+x_m$,  the evaluation map
$
 \ev_D:W\longrightarrow H^0(E|_D)
$
is an isomorphism.  We call such a system \emph{$m$-unisolvent}. We show that  on
$(\PP^1)^m$, the determinant of the universal ordered evaluation map has no
zeros away from the pairwise diagonals; moreover, it can be written as
$ \divisor(\delta_{E,W})
   =e\sum_{1\leq i<j\leq m}\Delta_{ij}
$ for some positive integer $e$.
We call $e$ the \emph{diagonal exponent}. The interesting part is that this single invariant $e$ controls three a priori different geometries: the
Pl\"ucker curve, the variety swept out by the moving projective spaces and
the degeneration of interpolation when points collide. In Section~\ref{sec:rank-two-applications}, we  prove the rank-two
classification.  As applications, we completely answer Day's question in a stronger form and prove a strengthened version of the critical Terracini problem.   Finally, in
Section~\ref{sec:unisolvent-category}, we show that the $m$-unisolvent
systems form an abelian category of finite length, that the Hermite
systems are simple and that the diagonal exponent is additive in short
exact sequences. We then study rank-three systems. For \(m=2\), the exponents \(3\) and \(5\) are completely controlled by the lower-rank factors, while \(7\) is the least exponent of a simple rank-three system.

We now describe the first group of results of the paper.

\begin{theoremA}
Let
$f:\PP^1\to\Gr(k,V)$ be an $m$-direct-sum curve of diagonal exponent $e$ where $m\text{ is at least }2$ and $\dim V\text{ is }mk$.
Then the following statements hold.
\begin{enumerate}[label=\textup{(\arabic*)}]
\item If the drill $d(f)$ is $1$, then $f$ is the \((k-1)\)-st
osculating curve of the
rational normal curve \(\nu_{km-1}\subset\PP^{km-1}\).

\item Let
$g:\PP^1\longrightarrow\PP(\textstyle\bigwedge^kV)$
be the Pl\"ucker map induced by $f$.
The curve $g(\PP^1)$ is a linearly normal rational normal curve of degree
$e(m-1)$ in its linear span.

\item  Write $U_t\text{ for the subspace corresponding to }f(t)$ and let
$X_f\text{ be the union}\bigcup_{t\in\PP^1}\PP(U_t)\subset\PP(V)$.
The natural morphism
$\nu:\PP_{\PP^1}(S)\longrightarrow X_f$
is the normalization of $X_f$.  Moreover, $X_f$ is an irreducible,
nondegenerate variety of dimension $k$ whose degree is $e(m-1)$.

\item
For the associated
$m$-unisolvent system $(E,W)$, let
$\mathcal H_m\text{ be the Hilbert scheme }\operatorname{Hilb}^m(\PP^1)\simeq\PP^m$ and let $E^{[m]}$ be the
length-$m$ tautological bundle of $E$ on $\mathcal H_m$.  Let
$\Delta_m\subset\mathcal H_m$ denote  the discriminant hypersurface of
nonreduced divisors.  The divisor $\divisor(\det\epsilon_{E,W})$ of
the universal evaluation morphism
$\epsilon_{E,W}:W\otimes\OO_{\mathcal H_m}\longrightarrow E^{[m]}$ is $\frac{e-k}{2}\,\Delta_m.$
\item
The exponent $e$ satisfies
 $k\leq e\leq
 \left\lfloor\frac{\binom{mk}{k}-1}{m-1}\right\rfloor,
  e\equiv k\pmod 2.$
Moreover, the following conditions are equivalent.
\begin{enumerate}[label=\textup{(\alph*)}]
\item $e=k$;
\item There is a direct-sum decomposition
$V=V_1\oplus\cdots\oplus V_k, \text{ with }\dim V_j=m$ for every $j$. There are
degree-$(m-1)$ rational normal embeddings
$c_j:\PP^1\to\PP(V_j)$. If $\ell_j(t)\subset V_j$ is the line
represented by $c_j(t)$, then
$U_t\text{ is }\ell_1(t)\oplus\cdots\oplus\ell_k(t)
\text{ for every }t\in\PP^1;$

\item The fiber of the universal evaluation morphism $\epsilon_{E,W}$
at $[D]$, namely $W\to H^0(E|_D)$, is an isomorphism for at least one
nonreduced closed subscheme $D\subset\PP^1$ of length $m$.

\end{enumerate}

\end{enumerate}
\end{theoremA}

The second group of results begins with the complete classification in
rank two.  This is the first genuinely Grassmannian case: the generic
Perkinson partition is either $(1,1)$ or $(2)$ and the two possibilities
lead to geometrically distinct models.

\begin{theoremB}
Let $m\geq2$.  The following statements hold.

 Every $m$-direct-sum curve
$f:\PP^1\to\Gr(2,V)$ with $\dim V=2m$ is, up to the action of
$\GL(V)$, exactly one of the following:
\begin{enumerate}[label=\textup{(\alph*)}]
\item there are complementary subspaces $V_1,V_2\subset V$, each of
dimension $m$ and degree-$(m-1)$ rational normal embeddings
$c_i:\PP^1\to\PP(V_i)$ such that
$U_t\text{ is }\ell_1(t)\oplus\ell_2(t)$ for every $t\in\PP^1$, where
$\ell_i(t)\subset V_i$ is the line represented by $c_i(t)$;
\item there is a degree-$(2m-1)$ rational normal embedding
$c:\PP^1\to\PP(V)$ with image $R$ such that
$\PP(U_t)\text{ is }T_{c(t)}R$ for every $t\in\PP^1$.
\end{enumerate}
The diagonal exponents of these two models are $2$ and $4$, respectively.

\end{theoremB}
 Theorem B above can be reformulated in a more geometrical way as follows.
\begin{theoremC}
    Let $m\geq2$, let $V$ be a complex vector space of dimension $2m$ and let
$X$ be a positive-dimensional  projective subvariety of $\Gr(2,V)$. Assume that $Y\subset \PP(V)$ is a closed subvariety  ruled  by lines over $X$. Suppose  $\langle\ell_{x_1},\ldots,\ell_{x_m}\rangle=\PP(V)$ for every
$m$ pairwise distinct closed points $x_1,\ldots,x_m\in X$ where $\ell_{x_i}$ denotes the line over $x_i$.
Then
$X$ is a smooth rational curve  and $Y$ is an irreducible nondegenerate surface.  Moreover, up to projective
equivalence, exactly one of the following cases occurs.
\begin{enumerate}[label=\textup{(\arabic*)},nosep]
\item  $Y$ is the balanced rational normal scroll $S_{m-1,m-1}$.  The isomorphism
$Y\simeq\PP^1\times\PP^1$ identifies $\OO_Y(1)$ with
$\OO_{\PP^1\times\PP^1}(m-1,1)$.  In particular, $Y$ is smooth
and $\deg Y\text{ is }2m-2$.
\item There is a degree-$(2m-1)$ rational normal embedding
$c:\PP^1\to\PP(V)$ with image $R$ such that
$\ell_t\text{ is }T_{c(t)}R$.  In this case
$Y\text{ is }\operatorname{Tan}(R)$ and its normalization is
$q:\PP^1\times\PP^1\to Y$ with
$q^*\OO_Y(1)\simeq\OO_{\PP^1\times\PP^1}(2m-2,1)$.  Moreover,
$\deg Y\text{ is }4m-4$ and  the variety $Y$ is nonnormal and its reduced singular locus is
$R$.
\end{enumerate}
Both families satisfy the stated spanning condition.
\end{theoremC}
 When $m$ is $2$, Theorem C gives a complete answer to Day's problem without  the generic hypothesis or any assumption on the dimension of $X$. Even when $\dim(X)$ is $1$, Theorem C is still essentially stronger, because we do not assume the smoothness of $X$ in advance.

 We also strengthen the critical Terracini theorem of
Ballico and Chiantini \cite[Proposition~7]{BallicoChiantini2024}.
 Let
$X\subset\PP(V)$ be a  smooth, projective and nondegenerate
variety of dimension $n$ $  (n\ge 1)$. We set $\rho_s(n):=\binom{n+s}{n}$ with  $s\geq1$. Let
$\widehat{\operatorname{Osc}}^s_xX\subset V$ denote the affine $s$-th
osculating space at $x$.  Assume that  $\dim V\text{ is }m\rho_s(n)$ for some integer $m\ge 2$. The reduced critical
$s$-th Terracini locus is
$\mathcal T_m^{(s),\mathrm{red}}(X):=\{(x_1,\ldots,x_m)\in
\operatorname{Conf}_m(X):\widehat{\operatorname{Osc}}^s_{x_1}X+
\cdots+\widehat{\operatorname{Osc}}^s_{x_m}X\neq V\}$, where
$\operatorname{Conf}_m(X)$ is the space of ordered $m$-tuples of pairwise
distinct points.

\begin{theoremD}
Suppose that $\dim V\text{ is }m\rho_s(n)$.  Then
$\mathcal T_m^{(s),\mathrm{red}}(X)$ is empty if and only if $n\text{ is }1$ and
$X\subset\PP(V)$ is projectively equivalent to the rational normal curve
$\nu_{m(s+1)-1}(\PP^1)\subset\PP^{m(s+1)-1}$.
\end{theoremD}

Using the same strategy as in  Theorem D, we prove a dual version of  critical Terracini rigidity.
Let $X\subset\PP^N=\PP(V)$ be a  smooth, projective and
nondegenerate subvariety of dimension $n\ge 1$. We set $c:=N-n$ and assume $N+1\text{ is }mc$ ($m\ge2$).

\begin{theoremE}
Suppose that $n\text{ is at least }1$ and the intersection satisfies
$T_{x_1}X\cap\cdots\cap T_{x_m}X=\varnothing$ for every $m$ pairwise
distinct points of $X$.  Then exactly one of the following occurs:
\begin{enumerate}[label=\textup{(\arabic*)},nosep]
\item $m$ is $2$, $N$ is $3$ and $X$ is the twisted cubic;
\item $m$ is $3$, $N$ is $2$ and $X$ is a smooth conic.
\end{enumerate}
Conversely, both curves satisfy the stated intersection condition.
\end{theoremE}

The final group of results describes how the low-rank classification
constrains rank-three systems. Here $L_m$ denotes the rank-one
Vandermonde system and $H_{m,k}$ the rank-$k$ Hermite system introduced
in Example~\ref{par:split-hermite-models}.

\begin{theoremF}
Fix $m\geq2$. The following statements hold.
\begin{enumerate}[label=\textup{(\arabic*)}]
\item The category $\mathcal U_m$ of $m$-unisolvent systems is abelian,
and every object has finite length. The Hermite system $H_{m,k}$ is
simple for every $k\geq1$.
\item A nonsimple rank-three object has diagonal exponent $3$ or $5$.
Every rank-three object of exponent $3$ is isomorphic to
$L_m^{\oplus3}$. Consequently, every rank-three object of exponent
greater than $5$ is simple.
\item For $m=2$, every rank-three object of exponent $5$ is isomorphic
to $H_{2,2}\oplus L_2$. The least diagonal exponent of a simple
rank-three $2$-unisolvent system is $7$.
\end{enumerate}
\end{theoremF}

\section{Vandermonde prototypes and \texorpdfstring{$m$}{m}-direct-sum curves}
\label{sec:affine-vandermonde}
\label{sec:direct-sum-varieties}

 Rational normal curves can be characterized by the following theorem. For the proof,  see \cite[Example~1.14 and Proposition~18.9]{HarrisAG} and
\cite[Theorem~1]{EisenbudHarrisMinimalDegree}.

\begin{theorem}\label{thm:projective-rnc-prototype}
Consider a nonconstant morphism
\(f:X\to\PP(V)\simeq\PP^{m-1}(m\geq2)\) from a smooth projective variety \(X\) and let \(\ell_x\subset V\) denote the line represented by  \(f(x)\). If
\(V\text{ decomposes as }\ell_{x_1}\oplus\cdots\oplus\ell_{x_m}\) for every ordered \(m\)-tuple of pairwise distinct
points \(x_1,\ldots,x_m\in X\), then \(X\text{ is isomorphic to }\PP^1\). Up to the action of
\(\GL(V)\), the morphism \(f\) is the rational normal embedding
\(\nu_{m-1}:\PP^1\to\PP^{m-1}\).
\end{theorem}

We next record the affine form of the same rigidity.  Its elementary proof is independent of the
projective characterization and has its own interest.

\begin{theorem}\label{thm:affine-vandermonde}
Suppose that \(m\text{ is at least }2\) and that \(u(t)=(p_0(t),\ldots,p_{m-1}(t))\) is a
polynomial row vector. For pairwise distinct complex numbers
\(t_1,\ldots,t_m\), the notation \(M_u(t_1,\ldots,t_m)\) denotes the matrix
whose \((i,j)\)-entry is \(p_j(t_i)\) \((1\leq i\leq m\),
\(0\leq j\leq m-1)\). If
\(M_u(t_1,\ldots,t_m)\) is nonsingular for every such tuple, then
\(p_0,\ldots,p_{m-1}\) span \(\CC[t]_{\leq m-1}\).  Equivalently, there is a matrix
\(B\in\GL_m(\CC)\) such that \(u(t)B\text{ is }(1,t,t^2,\ldots,t^{m-1})\).
\end{theorem}

\begin{proof}
With \(d=\max_j\deg p_j\), the degree-\(d\) homogenization of \(p_j\) in homogeneous coordinates
\([X:Y]\) on \(\PP^1\) is \(P_j(X,Y)=Y^d p_j(X/Y)\), where \(t\text{ is }X/Y\) on the open chart defined by
\(Y\neq0\).  Each \(P_j\) is a section of \(\OO_{\PP^1}(d)\) and their span is denoted by
\(W\subset H^0(\PP^1,\OO_{\PP^1}(d))\).

The polynomials \(p_0,\ldots,p_{m-1}\) are linearly independent;  otherwise, every evaluation
matrix would have linearly dependent columns.  Hence \(W\) has dimension \(m\).

The linear system \(W\) is base-point-free.  If all \(p_j\) vanished at a complex number \(a\), then
the row \(u(a)\) would be zero and every evaluation matrix containing \(a\) would be singular.  At
infinity, at least one polynomial has degree \(d\), so the corresponding homogeneous polynomial does
not vanish at \([1:0]\).  Thus \(W\) has no base point on \(\PP^1\).

Since \(W\subset H^0(\PP^1,\OO_{\PP^1}(d))\) has dimension \(m\), one has
\(m\leq d+1\).  For the reverse inequality, if
\(d\text{ is at least } m\), then a general section of the base-point-free linear system \(W\) on the smooth curve \(\PP^1\)
has a reduced zero divisor.  Moreover, infinity is not a base point, so the sections that do
not vanish there form a nonempty open subset of \(W\).  It follows that there is a section
\(s\in W\) whose zero divisor is reduced and supported on \(\Aff^1\).  This divisor consists of
\(d\) distinct points, any \(m\) of which may be denoted by \(t_1,\ldots,t_m\).

The nonzero section \(s\) then lies in the kernel of the evaluation map
\(W\to\bigoplus_{i=1}^m\OO_{\PP^1}(d)_{t_i}\).  In the basis
\(P_0,\ldots,P_{m-1}\), this map is represented by \(M_u(t_1,\ldots,t_m)\), which is nonsingular by
hypothesis.  This contradiction proves  \(d=m-1\).

Now \(W\) is an \(m\)-dimensional subspace of the \(m\)-dimensional vector space
\(H^0(\PP^1,\OO_{\PP^1}(m-1))\) and hence they are equal.  On the affine chart, the polynomials
\(p_0,\ldots,p_{m-1}\) therefore form a basis of \(\CC[t]_{\leq m-1}\).  Comparing this basis with
\(1,t,\ldots,t^{m-1}\) gives a matrix \(B\in\GL_m(\CC)\) such that
\(u(t)B=(1,t,\ldots,t^{m-1})\).
\end{proof}

We now introduce the Grassmannian analogue.

\begin{definition}
For a morphism \(f:X\to\Gr(k,V)\) from a projective variety to the Grassmannian of
\(k\)-planes in a vector space \(V\) of dimension \(km\), we let \(U_x\subset V\) denote the \(k\)-plane
represented by \(f(x)\).  The morphism \(f\) is an \emph{\(m\)-direct-sum morphism} if there is
\(V=U_{x_1}\oplus\cdots\oplus U_{x_m}\) for every ordered \(m\)-tuple of pairwise distinct points
\(x_1,\ldots,x_m\in X\). Throughout this article, unless stated otherwise,  $X$ is assumed to be \textbf{ positive-dimensional}. Equivalently, for every reduced $0-$cycle
\(x_1+\cdots+x_m\) of degree \(m\), the natural summation map
\(U_{x_1}\oplus\cdots\oplus U_{x_m}\to V\) is an isomorphism.
\end{definition}

It follows directly from the definition that  every \(m\)-direct-sum morphism is injective. There are also further restrictions on the source $X$.
\begin{proposition}\label{prop:dimension-bound}
Every \(m\)-direct-sum morphism \(f:X\to\Gr(k,V)\), where \(\dim V\text{ is }km\), satisfies
\(\dim X=1\).
\end{proposition}

\begin{proof}
For pairwise distinct points \(y_1,\ldots,y_{m-1}\) of \(X\), the subspace
\(W:=U_{y_1}+\cdots+U_{y_{m-1}}\) has dimension \((m-1)k\) by the
\(m\)-direct-sum property and hence has codimension \(k\) in \(V\).  The Schubert divisor
\(\Sigma_W\subset\Gr(k,V)\) consists of the \(k\)-planes that meet \(W\) nontrivially.  It is a
hyperplane section in the Pl\"ucker embedding.

If \(\dim X\text{ is at least }2\), \(f^{-1}(\Sigma_W)\) contains a positive-dimensional algebraic subset.
Hence
there is a point
\(y_m\notin\{y_1,\ldots,y_{m-1}\}\) satisfying \(f(y_m)\in\Sigma_W\).  Thus we have
\(U_{y_m}\cap W\neq\{0\}\). So \(U_{y_1}+\cdots+U_{y_m}\) has dimension strictly smaller than \(km\).
This contradicts the \(m\)-direct-sum property.
\end{proof}

\begin{proposition}\label{prop:curve-is-p1}
If a smooth projective curve \(X\) admits an  \(m\)-direct-sum morphism
\(f:X\to\Gr(k,V)\), then \(X\text{ is isomorphic to }\PP^1\).
\end{proposition}

\begin{proof}
The pullback of the tautological rank-\(k\) subbundle on \(\Gr(k,V)\) is denoted by
\(S=f^*\mathcal U\) and we let \(E\text{ be }S^\vee\).  The fiber \(S_x\) is \(U_x\subset V\).  Dualizing the
tautological inclusion gives a quotient \(V^\vee\otimes\OO_X\to E\). On \(X^m\), the determinant of
the evaluation morphism \(V^\vee\otimes\OO_{X^m}\to\bigoplus_{i=1}^m p_i^*E\) is a section
\(s\in H^0(X^m,\bigotimes_{i=1}^m p_i^*\det E)\).

At \((x_1,\ldots,x_m)\in X^m\), the fiber map is
\(V^\vee\to U_{x_1}^\vee\oplus\cdots\oplus U_{x_m}^\vee\), dual to the summation map
\(U_{x_1}\oplus\cdots\oplus U_{x_m}\to V\).  Hence \(s(x_1,\ldots,x_m)\)  is not zero precisely when
this summation map is an isomorphism.  The \(m\)-direct-sum property implies that this holds whenever
the \(x_i\) are pairwise distinct.  Thus the support of \(\divisor(s)\) is contained in the union of
the pairwise diagonals \(\Delta_{ij}\subset X^m\).

Since \(X\) is a smooth irreducible curve, \(X^m\) is smooth and every \(\Delta_{ij}\) is an
irreducible Cartier divisor.  There are therefore nonnegative integers \(c_{ij}\) satisfying
\(\divisor(s)=\sum_{1\leq i<j\leq m}c_{ij}\Delta_{ij}\).

This divisor is symmetric.  A permutation of the factors of \(X^m\) permutes the rank-\(k\) summands
of the target of the evaluation map.  Although the determinant section may be multiplied by the sign
arising from the permutation of blocks of size \(k\), its divisor is unchanged.  The group \(S_m\)
acts transitively on the unordered pairs \(\{i,j\}\), so all \(c_{ij}\) are equal.  Consequently,
the divisor has the form \(\divisor(s)=c\sum_{1\leq i<j\leq m}\Delta_{ij}\) for an integer
\(c\geq0\).

The integer \(c\) is positive.  Since \(f\) is nonconstant, the Pl\"ucker line bundle
\(\det E=f^*\OO_{\Gr(k,V)}(1)\) has positive degree.  For distinct points
\(y_2,\ldots,y_m\) of \(X\), the restriction of \(s\) to the slice
\(X\times\{y_2\}\times\cdots\times\{y_m\}\) is a section of \(\det E\) with zero divisor
\(c(y_2+\cdots+y_m)\).  Hence
\(\det E\text{ is isomorphic to }\OO_X(c(y_2+\cdots+y_m))\) and the positivity of \(\deg\det E\) implies  \(c>0\).

For fixed distinct points \(y_2,\ldots,y_{m-1}\) and any points
\(y,y'\in X\setminus\{y_2,\ldots,y_{m-1}\}\), we have the following isomorphisms
\(\det E\simeq\OO_X(c(y_2+\cdots+y_{m-1}+y))
\simeq\OO_X(c(y_2+\cdots+y_{m-1}+y'))\).  It follows that
\(\OO_X(cy)\text{ is isomorphic to }\OO_X(cy')\), which is equivalent to \(cy\sim cy'\).

If \(g(X)\) is positive, the Abel--Jacobi morphism based at a fixed point \(x_0\in X\) is
\(\alpha:X\to\Jac(X)\),
\(x\mapsto[\OO_X(x-x_0)]\).  For a smooth projective curve of positive genus, this morphism is
nonconstant.  Multiplication by \(c\) on \(\Jac(X)\) is a finite morphism, so the composite
\([c]\circ\alpha\) is also nonconstant and sends \(x\) to
\([\OO_X(c(x-x_0))]\).  The relation \(cy\sim cy'\) on a dense open subset instead makes this
composite constant on a dense open subset, and hence on all of \(X\). Thus we get  a contradiction.  Therefore
\(X\) is isomorphic to \(\PP^1\).
\end{proof}
\begin{theorem}\label{source}
  The source of every $m$-direct-sum morphism is a smooth rational curve. Thus the smoothness  assumption on $X$ in Proposition \ref{prop:curve-is-p1} can be removed.\end{theorem}
  \begin{proof}
      We use the notation from the proof of Proposition \ref{prop:curve-is-p1}.
  Let
$\nu:C\to X$ be the normalization, set $T=\nu^*f^*\mathcal U$, and write
$g=f\circ\nu$.  On $C^m$, with projections $p_i:C^m\to C$, let $\delta$
be the determinant of the bundle map
$\bigoplus_{i=1}^m p_i^*T\to V\otimes\OO_{C^m}$.  This is a nonzero section
of a line bundle.  Let $W\subset C^m$ be the complement of the pairwise
diagonals.  The direct-sum hypothesis implies that the zero locus of
$\delta|_W$ is contained in the union, over $i<j$, of the inverse images of
the off-diagonal part of $C\times_XC$ under $(p_i,p_j)$.  This off-diagonal
part is finite. The indicated union therefore has
codimension two in $W$.  On the other hand, $W$ is smooth and integral, so
the nonempty zero locus of a nonzero section of a line bundle on $W$ would
be an effective Cartier divisor and hence have pure codimension one.  Thus
$\delta|_W$ is nowhere zero and $g$ is itself an $m$-direct-sum morphism. By
Proposition~\ref{prop:curve-is-p1}, $C\text{ is isomorphic to }\PP^1$. Then, by Theorem \ref{thm:plucker-vandermonde}  below, the composition of the following morphisms $C\overset{\nu} \to X \overset{f}\to Gr(k,V)\xrightarrow[]{\text{Pl\"ucker}} \PP(\bigwedge^k V)$ is a closed embedding. So $\nu$ is an isomorphism and hence $X$ is smooth.
\end{proof}
\section{Perkinson normal forms and extremal classifications}\label{sec:perkinson}

The preceding source-rigidity results reduce the higher-dimensional problem to curves.  Before introducing the diagonal exponent, we recall the rough
differential classification inspired by the derived quotients in Perkinson's study
of algebraic curves in Grassmannians \cite{Perkinson}.

\subsection{Derived quotients and the drill}

\textbf{Notation:} We first fix the notation used throughout this paper. For a nonconstant morphism \(f:\PP^1\to\Gr(k,V)\),  let \(U_x\subset V\) denote the \(k\)-plane
represented by \(f(x)\).  Let $\mathcal U$ be the universal subbundle on  $\Gr(k,V)$  and let \(S=f^*\mathcal U\) be the pullback of $\mathcal{U}$.  Assume
\(E \text{ is } S^\vee \), then there is a quotient $V^\vee\otimes \mathcal{O}_{\PP^1}\to E$.

If
\(S_{E}\) denotes the kernel of the quotient
\(V^\vee\otimes\OO_{\PP^1}\to E\), then the tangent map gives a morphism
\(\partial:S_{E}\to\Omega^1_{\PP^1}\otimes E\).  The sequence
\(0\to S_{E}\to V^\vee\otimes\OO_{\PP^1}\to E\to0\) is exact and Perkinson's first derived
quotient is \(E_1=(\coker(\partial)\otimes(\Omega^1_{\PP^1})^{-1})/\tors\).

The composition \(V^\vee\otimes\OO_{\PP^1}\to E\to E_1\) is again a quotient of the trivial
bundle and therefore determines a new Grassmannian curve.  Geometrically, this operation removes the
first-order moving directions of the original curve. For consistency of labeling, we let $E_0$ be $E$.  Iteration gives the sequence of derived quotients
\(E_0\longrightarrow E_1\longrightarrow E_2\longrightarrow\cdots\).

The \(i\)-th differential rank is
\(\drk_i(f)=\rk E_{i-1}-\rk E_i\) and the \emph{drill} is \(d(f)=\drk_1(f)\).  The terminal quotient
is a quotient of \(V^\vee\otimes\OO_{\PP^1}\) whose dual is the fixed space
\(L(f)=\bigcap_{x\in\PP^1}U_x\subset V\).

Passage to the quotient by \(L(f)\) gives the quotient  \(\overline U_t=U_t/L(f)\).  On an affine chart of \(\PP^1\),
Perkinson's normal form provides vector-valued functions \(u_1(t),\ldots,u_d(t)\) and nonnegative
integers \(a_1,\ldots,a_d\) such that \(\overline U_t\) is the direct sum, over
\(1\leq\alpha\leq d\), of the spaces
\(\Span\{u_\alpha(t),u_\alpha'(t),\ldots,u_\alpha^{(a_\alpha)}(t)\}\).
These integers satisfy \(k-\dim L(f)=\sum_{\alpha=1}^d(a_\alpha+1)\).

\begin{remark}\label{prop:perkinson}
Following  the preceding notation, Perkinson's normal-form theorem
\cite[Theorem~5.1 and Remark~6.1.2]{Perkinson} identifies the number of strings \(d\) with the drill
\(d(f)\).  Moreover, the differential ranks are
\(\drk_i(f)=\#\{\alpha:a_\alpha\geq i-1\}\) and in the derived sequence above, the dual of \(E_i\) corresponds to the
subbundle generated by the direct sum, over \(1\leq\alpha\leq d\), of the spaces
\(\Span\{u_\alpha(t),u_\alpha'(t),\ldots,u_\alpha^{(a_\alpha-i)}(t)\}\),
where the summands with \(a_\alpha<i\) are omitted.
\end{remark}

For an \(m\)-direct-sum curve in \(\Gr(k,km)\), the fixed space \(L(f)\) is zero: a nonzero vector in
that intersection would belong to every fiber and contradict the directness of any two fibers. This
shows that the integers \(n_\alpha=a_\alpha+1\) form a partition of \(k\).

\subsection{The case \texorpdfstring{$d(f)=1$}{d(f)=1}}

\begin{theorem}\label{thm:single-string}
Every \(m\)-direct-sum curve \(f:\PP^1\to\Gr(k,V)\) with \(\dim V=km\)
and \(d(f)=1\) is, up to the action of \(\GL(V)\), the \((k-1)\)-st
osculating curve of the
rational normal curve \(\nu_{km-1}\subset\PP^{km-1}\).
\end{theorem}
\begin{proof}
We write $N:=km$. Since $d(f)$ is $1$, Perkinson's normal form has a single
string. Since the fixed space is zero, as noted in Remark~2.1,
[9, Corollary~6.2.2] shows that the last nonzero derived quotient is a line
quotient $q:V^\vee\otimes\mathcal O_{\mathbb P^1}\to L$. After dualization,
$L^\vee$ is a line subbundle of $V\otimes\mathcal O_{\mathbb P^1}$ and therefore
defines a morphism $\varphi:\mathbb P^1\to\mathbb P(V)$ with
$L=\varphi^*\mathcal O_{\mathbb P(V)}(1)$. On a local coordinate chart with
coordinate $t$, a generator of $L^\vee$ is denoted by $v(t)$.
By [9, Theorem~5.1], there is a dense open subset
$U^\circ\subset\mathbb P^1$ on which Perkinson's normal form gives
$U_t=\operatorname{Span}\{v(t),v'(t),\ldots,v^{(k-1)}(t)\}$.
Thus $f$ is locally the $(k-1)$-st osculating curve of $\varphi$ on $U^\circ$.
We do not yet assume that the ordinary derivatives form a frame at the
remaining points.

We prove that $\varphi$ is the rational normal curve of degree $N-1$.
Every $\lambda\in V^\vee$ determines a global section $q(\lambda)$ of $L$
and its $(N-1)$-jet defines a bundle morphism
$\Theta:V^\vee\otimes\mathcal O_{\mathbb P^1}\to J^{N-1}L$.
Both bundles have rank $N$, so $\Theta$ is an isomorphism if and only if
its determinant is nowhere zero.

We first describe this determinant locally. On a sufficiently small open
subset $A\subset\mathbb P^1$ with coordinate $t$, the line bundle $L$ has
a local frame $e$ and one may choose a local lift $v(t)\in V$ satisfying
$q(\lambda)=\lambda(v(t))e$ for every $\lambda\in V^\vee$.
In this coordinate and frame, the $(N-1)$-jet of a section $g(t)e$ is
represented by $(g(t),g'(t),\ldots,g^{(N-1)}(t))$.
Hence $\Theta$ sends $\lambda$ locally to
$(\lambda(v(t)),\lambda(v'(t)),\ldots,\lambda(v^{(N-1)}(t)))$.

For a basis $\lambda_1,\ldots,\lambda_N$ of $V^\vee$ and its dual basis
$e_1,\ldots,e_N$ in $V$, the local matrix of $\Theta$ has entries
$M_{a,j}(t)=\lambda_j(v^{(a)}(t))$ for $0\leq a\leq N-1$.
If $v^{(a)}(t)$ is $\sum_{j=1}^N x_{j,a}(t)e_j$, then $x_{j,a}(t)$ is
$\lambda_j(v^{(a)}(t))$. Thus $M(t)$ is the transpose of the matrix whose
columns are $v(t),v'(t),\ldots,v^{(N-1)}(t)$.
Up to a nonzero constant determined by the basis and the convention for
jet coordinates, $\det\Theta$ is locally represented by the Wronskian
$\mathcal W(t)=\det(v(t),v'(t),\ldots,v^{(N-1)}(t))$.

The function $\mathcal W$ is not globally defined; it is a local
representative of the global section $\det\Theta$.
Replacing $v(t)$ by $h(t)v(t)$ multiplies $\mathcal W(t)$ by $h(t)^N$
and changing the coordinate from $t$ to $\tau$ multiplies it by a power
of $dt/d\tau$. These factors are units, so the zero loci of the local
Wronskians glue to the zero divisor of $\det\Theta$ on $U_0$.

The Wronskian is not identically zero. Otherwise, after writing the
coordinates of $v(t)$ in a basis of $V$ as $f_1(t),\ldots,f_N(t)$, their
ordinary Wronskian would vanish identically. In characteristic zero, the
Wronskian criterion would give a nonzero $\lambda\in V^\vee$ such that
$\lambda(v(t))=0$ throughout the chart. Differentiation would then give
$\lambda(v'(t))=\cdots=\lambda(v^{(k-1)}(t))=0$ and hence
$U_t\subset\ker\lambda$ on the intersection of the chart with
$U^\circ$. The sum of the fibers over any $m$ distinct points
in this intersection would lie in $\ker\lambda$, contradicting
the $m$-direct-sum property.

Let $\delta$ be the intrinsic $m$-point determinant of the
summation map
$\bigoplus_{i=1}^m p_i^*S\to V\otimes\mathcal O_{(\mathbb P^1)^m}$, where
$S=f^*\mathcal U$ and the $p_i$ are the projections. As in the proof of
Proposition~1.5, the direct-sum property and permutation symmetry give
$\operatorname{div}(\delta)=c\sum_{i<j}\Delta_{ij}$ for some $c\geq1$.
Local representatives of $\delta$ are computed using genuine frames of
$S$.

At a point $t_0$ in $U^\circ$ for which $\mathcal W(t_0)\neq0$,
the local $m$-point determinant $D(t_1,\ldots,t_m)$ is defined as
$\det\bigl(v(t_i),v'(t_i),\ldots,v^{(k-1)}(t_i)\bigr)_{i=1}^m$
and it is nonzero whenever the $t_i$ are pairwise distinct.
Here the coordinate neighborhood is chosen inside $U^\circ$,
so $D$ is a local representative of $\delta$.
Its zero divisor is therefore $c\sum_{i<j}\Delta_{ij}$.
We determine $c$ at $t_0$. Writing $t_i=t_0+\varepsilon_i$, Taylor
expansion for $0\leq a\leq k-1$ gives
$v^{(a)}(t_0+\varepsilon_i)=\sum_{b\geq0}\varepsilon_i^b
v^{(a+b)}(t_0)/b!$.

It follows that the leading term of
$D(t_0+\varepsilon_1,\ldots,t_0+\varepsilon_m)$ is
$C\prod_{i<j}(\varepsilon_i-\varepsilon_j)^{k^2}\mathcal W(t_0)$,
while all remaining terms have higher order; here $C\neq0$ depends only
on $k$ and $m$.
For completeness, a nonzero determinant term must use $N$
distinct derivative orders. The smallest possible total degree in the
$\varepsilon_i$ is therefore
$
\ell=\frac{N(N-1)}2-m\frac{k(k-1)}2
     =k^2\binom m2.
$
The homogeneous part of degree $\ell$ is
$\mathcal W(t_0)\det Q(\varepsilon)$, where the rows of $Q$ are indexed by
$0\leq q<N$, the columns by $(i,a)$ with $1\leq i\leq m$ and $0\leq a<k$,
and
$
Q_{q,(i,a)}(\varepsilon)=
\begin{cases}
\varepsilon_i^{q-a}/(q-a)!,&q\geq a,\\
0,&q<a.
\end{cases}
$
For pairwise distinct $\varepsilon_i$, the transpose of $Q$ is the
Hermite evaluation matrix on polynomials of degree at most $N-1$, in
the basis $x^q/q!$. Its kernel is zero: a polynomial in the kernel
would be divisible by $\prod_i(x-\varepsilon_i)^k$, of degree $N$.
Thus $\det Q$ vanishes only on the diagonals. Its divisor is symmetric,
and it is homogeneous of degree $\ell$; hence
$\det Q=C\prod_{i<j}(\varepsilon_i-\varepsilon_j)^{k^2}$, as asserted.
Since $\mathcal W(t_0)\neq0$, the leading term is nonzero, and hence
$c=k^2$. Thus the global intrinsic determinant satisfies
$
\operatorname{div}(\delta)=k^2\sum_{i<j}\Delta_{ij}.
$

With $\Delta(t_1,\ldots,t_m)=\prod_{i<j}(t_i-t_j)$, the quotient
$F:=D/\Delta^{k^2}$ extends holomorphically across the diagonals and has
no zeros, so it is a unit on the local product chart
contained in $(U^\circ)^m$.
We use the same notation for its restriction to the small diagonal.
Writing $t_i=t+\varepsilon_i$, the same Taylor calculation shows that
the homogeneous term of degree $\ell$ of
$D(t+\varepsilon_1,\ldots,t+\varepsilon_m)$ is
$C\prod_{i<j}(\varepsilon_i-\varepsilon_j)^{k^2}\mathcal W(t)$,
while all remaining terms have higher order.
Since $D=\Delta^{k^2}F$, comparison of the lowest homogeneous terms in
the variables $\varepsilon_i$ yields
$F(t,\ldots,t)=C\mathcal W(t)$.
The function $F$ is nowhere zero, so $\mathcal W(t)$ is nowhere zero on
the chart.

We now treat an arbitrary point $p\in\mathbb P^1$, including
those outside $U^\circ$. Choose a coordinate $z$ centered at $p$, a local
lift $v(z)$ as above, and a genuine local frame
$B(z)=(b_1(z),\ldots,b_k(z))$ of $S$.
The vectors $v(z),v'(z),\ldots,v^{(k-1)}(z)$ are sections of $S$:
their images in $(V\otimes\mathcal O)/S$ vanish on the dense open subset
$U^\circ$, hence vanish identically. Consequently,
$
R(z):=(v(z),v'(z),\ldots,v^{(k-1)}(z))=B(z)A(z)
$
for a holomorphic $k\times k$ matrix $A(z)$.
Put $h(z):=\det A(z)$; it is not identically zero, but it may a priori
vanish at $p$.

The determinant formed from the genuine frames is
$
\widetilde D(z_1,\ldots,z_m)
=\det(B(z_1),\ldots,B(z_m))
=\Delta(z_1,\ldots,z_m)^{k^2}H(z_1,\ldots,z_m),
$
where $H$ is a unit by the divisor formula for $\delta$.
The ordinary-derivative determinant therefore satisfies
$
D(z_1,\ldots,z_m)
=\left(\prod_{i=1}^m h(z_i)\right)\widetilde D(z_1,\ldots,z_m).
$
In particular, $F=D/\Delta^{k^2}$ is holomorphic on this product chart,
although it has not yet been shown to be a unit.
The same homogeneous-term calculation, valid even where
$\mathcal W(z)=0$, gives
$
C\mathcal W(z)=F(z,\ldots,z)=h(z)^m H(z,\ldots,z).
$
Thus, with $\beta:=\operatorname{ord}_p\mathcal W$ and
$\gamma:=\operatorname{ord}_p h$, we obtain
$
\beta=m\gamma.
$

We show that this equality forces $\beta=\gamma=0$.
After a constant change of basis of $V$, write the coordinate functions
of $v(z)$ as $f_0(z),\ldots,f_{N-1}(z)$, with
$
f_j(z)=u_j z^{a_j}+O(z^{a_j+1}),
u_j\neq0,  0=a_0<a_1<\cdots<a_{N-1}.
$
Such a basis exists by elimination of leading Taylor coefficients,
since the coordinate functions are linearly independent.
Moreover, $a_0=0$ because $q$ is surjective and hence $v(p)\neq0$.
Set $\alpha_j:=a_j-j$. Then
$0=\alpha_0\leq\alpha_1\leq\cdots\leq\alpha_{N-1}$.
For any $1\leq r\leq N$ and $j_0<\cdots<j_{r-1}$,
$
\operatorname{ord}_p\det\bigl(f_{j_\nu}^{(a)}(z)\bigr)_{0\leq\nu,a<r}
=\sum_{\nu=0}^{r-1}a_{j_\nu}-\frac{r(r-1)}2.
$
Indeed, the coefficient of this power of $z$ is
$
\left(\prod_{\nu=0}^{r-1}u_{j_\nu}\right)
\det\bigl((a_{j_\nu})_a\bigr)_{0\leq\nu,a<r}
=\left(\prod_{\nu=0}^{r-1}u_{j_\nu}\right)
\prod_{\mu<\nu}(a_{j_\nu}-a_{j_\mu})\neq0,
$
where $(x)_a=x(x-1)\cdots(x-a+1)$ and $(x)_0=1$.
Taking $r=N$ gives
$\beta=\sum_{j=0}^{N-1}\alpha_j$.
For $r=k$, each maximal minor of $R=BA$ is $h$ times the corresponding
minor of $B$, and at least one maximal minor of $B$ is a unit at $p$.
Thus $\gamma$ is the smallest order of a maximal minor of $R$, attained
by the rows $0,\ldots,k-1$, and
$\gamma=\sum_{j=0}^{k-1}\alpha_j$.

Since $N=mk$, the equality $\beta=m\gamma$ now reads
$
0=\sum_{q=1}^{m-1}\sum_{j=0}^{k-1}
       (\alpha_{qk+j}-\alpha_j).
$
Every summand is nonnegative. As $m\geq2$, the summand with $q=1$ and
$j=0$ shows that $\alpha_k=\alpha_0=0$.
Monotonicity then gives $\alpha_0=\cdots=\alpha_{k-1}=0$;
hence $\gamma=0$ and $\beta=m\gamma=0$.
Therefore $h(p)\neq0$ and $\mathcal W(p)\neq0$.
Since $p$ was arbitrary, this also treats every point omitted from
$U^\circ$.

It follows globally that $\det\Theta$ has no zeros and that
$\Theta:V^\vee\otimes\mathcal O_{\mathbb P^1}\to J^{N-1}L$ is an
isomorphism. In particular, $J^{N-1}L$ is trivial.

The line bundle $L$ has the form $L\simeq\mathcal O_{\mathbb P^1}(b)$.
The standard filtration of the jet bundle gives
$\deg J^{N-1}L=\sum_{a=0}^{N-1}(b-2a)$.
This sum equals $Nb-N(N-1)=N\bigl(b-(N-1)\bigr)$.
Since $J^{N-1}L$ is trivial, we have identities $b=N-1=km-1$ and
therefore $L$ is isomorphic to $\mathcal O_{\mathbb P^1}(km-1)$.

Finally, the quotient $V^\vee\otimes\mathcal O_{\mathbb P^1}\to L$
induces an injective map $V^\vee\to H^0(\mathbb P^1,L)$.
Indeed, if a nonzero $\lambda\in V^\vee$ mapped to the zero section,
then every local lift $v(t)$ satisfies $\lambda(v(t))=0$.
And hence $U_t\subset\ker\lambda$ for all $t$, which is impossible.
Both $V^\vee$ and $H^0(\mathbb P^1,\mathcal O_{\mathbb P^1}(km-1))$
have dimension $km$, so this injection is an isomorphism.
Thus $\varphi$ is defined by the complete linear system
$|\mathcal O_{\mathbb P^1}(km-1)|$ and is the rational normal curve of
degree $km-1$.
Since $h(p)\neq0$ at every point, the ordinary derivatives
now form a frame of $S$ everywhere.
Perkinson's local normal form now identifies $f$ with its $(k-1)$-st
osculating curve.
\end{proof}
 \section{The diagonal exponent and its geometry}\label{sec:unisolvent}

 This section recasts the direct-sum
condition as an intrinsic interpolation problem.  The resulting determinant defines the
diagonal exponent \(e\), which refines the differential classification and governs three
features of the geometry: its bounds and extremal behavior, the degree of the variety
swept out by the moving subspaces and the failure of interpolation when reduced
divisors collide.

\begin{definition}\label{def:coherent-system}
A coherent system on \(\PP^1\) is a pair \((E,W)\), where \(E\) is a vector bundle on
\(\PP^1\) and \(W\subset H^0(\PP^1,E)\) is a finite-dimensional subspace.  For a
reduced divisor \(D=x_1+\cdots+x_m\), the evaluation map is
\(\ev_D:W\to H^0(E|_D)=\bigoplus_{i=1}^m E_{x_i}\).
\end{definition}

\begin{definition}
A coherent system \((E,W)\) is \(m\)-unisolvent of rank \(k\) if \(\rk E\text{ is }k\),
\(\dim W\text{ is }mk\) and \(\ev_D\) is an isomorphism for every reduced divisor \(D\) of
degree \(m\).
\end{definition}

\begin{proposition}\label{prop:equivalence}
There is a natural bijection between isomorphism classes of rank-\(k\)
\(m\)-unisolvent systems \((E,W)\) on \(\PP^1\) and \(m\)-direct-sum curves
\(f:\PP^1\to\Gr(k,W^\vee)\), up to the action of \(\GL(W)\).  Under this
bijection, \(f(x)\) is the image of the inclusion \(E_x^\vee\hookrightarrow W^\vee\).
\end{proposition}

\begin{proof}
For a rank-\(k\) \(m\)-unisolvent system \((E,W)\) and any \(x\in\PP^1\), there exist
points \(x_2,\ldots,x_m\) distinct from one another and from \(x\).  Since
\(\ev_{x+x_2+\cdots+x_m}\) is an isomorphism, the map \(W\to E_x\) is surjective.
It follows that \(W\otimes\OO_{\PP^1}\twoheadrightarrow E\) is surjective.  The dual
inclusion \(E^\vee\hookrightarrow W^\vee\otimes\OO_{\PP^1}\) determines a curve in
\(\Gr(k,W^\vee)\) and the evaluation isomorphism for \(D\) is precisely the
direct-sum condition for the corresponding \(m\) subspaces.

Conversely, for an \(m\)-direct-sum curve \(f:\PP^1\to\Gr(k,V)\), the dual of the
pullback of the tautological subbundle is denoted by \(E\).  Dualizing the tautological
inclusion gives a surjection \(V^\vee\otimes\OO_{\PP^1}\twoheadrightarrow E\).  The
induced map \(V^\vee\to H^0(E)\) is injective: a nonzero kernel element would vanish on
every subspace \(f(x)\), whereas any \(m\) fibers over distinct points span \(V\).
Taking \(W=V^\vee\) therefore gives an \(m\)-unisolvent system.  The two constructions
are inverse and changes of basis in \(W\) correspond to the action of \(\GL(W)\).
\end{proof}

\subsection{The determinant criterion and the exponent}

For a coherent system \((E,W)\) with \(\rk E=k\) and \(\dim W=km\), a choice of a
nonzero element of \(\det W\) turns the determinant of the evaluation morphism
\(W\otimes\OO_{(\PP^1)^m}\to\bigoplus_{i=1}^m p_i^*E\) into a section
\(\delta_{E,W}\in H^0((\PP^1)^m,\bigotimes_{i=1}^m p_i^*\det E)\).  The diagonal defined by
\(x_i=x_j\) is denoted by \(\Delta_{ij}\subset(\PP^1)^m\).

\begin{proposition}\label{prop:detcriterion}
A coherent system \((E,W)\) is \(m\)-unisolvent if and only if
its determinant divisor is \(e\sum_{1\leq i<j\leq m}\Delta_{ij}\) for some
integer \(e\geq1\).  In this case,  we have
\(\deg E=e(m-1)\) and \(e\equiv k^2\pmod 2\).  The integer \(e\) is called the
\emph{diagonal exponent} of \((E,W)\).
\end{proposition}

\begin{proof}
The \(m\)-unisolvent condition says exactly that the determinant section vanishes only
when two of the points coincide.  Symmetry gives the stated divisor.  Restricting it
to \(\PP^1\times\{x_2\}\times\cdots\times\{x_m\}\) gives
\(\det E\simeq\OO_{\PP^1}(e(x_2+\cdots+x_m))\) and hence
\(\deg E\text{ is }e(m-1)\).  Transposing two blocks of size \(k\) changes the determinant by
\((-1)^{k^2}\), while \(\prod_{i<j}(t_i-t_j)^e\) changes by \((-1)^e\).  This proves
the parity assertion.
\end{proof}

The interpolation condition also imposes a useful positivity constraint on the
underlying bundle.

\begin{lemma}\label{lem:unisolvent-splitting-lower-bound}
Let  \((E,W)\) be an \(m\)-unisolvent system of rank \(r\). If
\(E\text{ decomposes as }\bigoplus_{i=1}^r\OO_{\PP^1}(a_i)\), then \(a_i\text{ is at least } m-1\) for every
\(i\).  Consequently, both \(H^1(E)\) and \(H^1(E(-m))\) vanish.
\end{lemma}

\begin{proof}
Unisolvence implies that \(W\otimes\OO_{\PP^1}\to E\) is surjective, so \(E\) is
globally generated. Hence every \(a_i\) is nonnegative.  For a reduced divisor \(D\)
of degree \(m\), the map \(W\to H^0(E|_D)\) is an isomorphism, and hence the full
restriction map \(H^0(E)\to H^0(E|_D)\) is surjective.  The cohomology sequence of
\(0\to E(-D)\to E\to E|_D\to0\), together with \(H^1(E)=0\), gives
\(H^1(E(-D))=0\).  Since \(E(-D)\simeq E(-m)\),  we have
\(H^1(\OO_{\PP^1}(a_i-m))=0\) for every \(i\), which is equivalent to
\(a_i\geq m-1\).
\end{proof}

For an \(m\)-direct-sum curve \(f:\PP^1\to\Gr(k,V)\), we let the Pl\"ucker embedding  induced by $f$ be
\(g:\PP^1\to\PP(\bigwedge^k V)\). Let the linear span of $g$ be
\(L_f=\langle g(\PP^1)\rangle\).

\begin{theorem}\label{thm:plucker-vandermonde}
If \(e\) is the diagonal exponent of \(f\) and if \(d=e(m-1)\), then
\(g(\PP^1)\subset L_f\simeq\PP^d\) is a linearly normal rational normal curve of
degree \(d\).
\end{theorem}

\begin{proof}
The bundles \(S=f^*\mathcal U\) and \(E=S^\vee\) satisfy
\(\det E\simeq\OO_{\PP^1}(d)\). The Pl\"ucker coordinates define a linear map
\(\rho:\bigwedge^k V^\vee\to H^0(\PP^1,\det E)=H^0(\PP^1,\OO_{\PP^1}(d))\).
We now prove that \(\rho\) is surjective.

For a reduced divisor \(D=y_1+\cdots+y_{m-1}\), the direct-sum condition makes
\(A_D=U_{y_1}\oplus\cdots\oplus U_{y_{m-1}}\) a \(k(m-1)\)-dimensional
subspace of \(V\).  A volume form \(\omega\in\bigwedge^{mk}V^\vee\) and nonzero
Pl\"ucker vectors \(p_{y_i}\in\bigwedge^k U_{y_i}\) determine the Pl\"ucker linear form
\(\ell_D\) given by
\(\ell_D(z)=\omega(p_{y_1}\wedge\cdots\wedge p_{y_{m-1}}\wedge z)\) for
\(z\in\bigwedge^k V\).
Its pullback by \(g\) is the \(m\)-point evaluation determinant with the last
\(m-1\) entries fixed.  Proposition~\ref{prop:detcriterion} therefore gives
\(\divisor(\rho(\ell_D))=e(y_1+\cdots+y_{m-1})=eD\).  Consequently, if
\(q_D\in H^0(\OO_{\PP^1}(m-1))\) has zero divisor \(D\), then \(\rho(\ell_D)\) is
proportional to \(q_D^e\).

The pure powers \(q^e\), with \(q\in H^0(\OO_{\PP^1}(m-1))\), span
\(H^0(\OO_{\PP^1}(e(m-1)))\).  Indeed, polarization expresses every product
\(q_1\cdots q_e\) as a linear combination of pure powers.  Such products span all
binary forms of degree \(e(m-1)\): every binary monomial of that degree is a product
of \(e\) binary monomials of degree \(m-1\), because every integer between \(0\) and
\(e(m-1)\) is a sum of \(e\) integers between \(0\) and \(m-1\).  Reduced divisors
form a dense open subset of \(\PP H^0(\OO_{\PP^1}(m-1))\), so the powers \(q_D^e\)
with \(D\) reduced already have the same span.  Hence \(\rho\) is surjective.

The Pl\"ucker map is consequently defined by the complete linear system
\(|\OO_{\PP^1}(d)|\).  It is the degree-\(d\) rational normal embedding and its
linear span has dimension \(d\).
\end{proof}

The exponent also records the differential type of the curve.

\begin{proposition}\label{prop:lower-bound}
For a rank-\(k\) \(m\)-unisolvent system \((E,W)\), if the Perkinson normal form of
the corresponding \(m\)-direct-sum curve has partition \((n_1,\ldots,n_s)\) at a
general point of \(\PP^1\), where \(\sum_{\alpha=1}^s n_\alpha=k\), then its
diagonal exponent satisfies
\(e\geq\sum_{\alpha=1}^s n_\alpha^2\).
\end{proposition}

\begin{proof}
For the corresponding \(m\)-direct-sum curve
\(f:\PP^1\to\Gr(k,W^\vee)\), we write \(U_t:=f(t)\).  At a general point \(p\) where
the Perkinson normal form is valid, a local coordinate \(z\) centered at \(p\) may be chosen.
On a sufficiently small coordinate neighborhood, the normal form gives a holomorphic frame of \(U_z\) consisting of
\(u_\alpha(z),u_\alpha'(z),\ldots,u_\alpha^{(n_\alpha-1)}(z)\) for
\(1\leq\alpha\leq s\).

Pairwise distinct points \(q_3,\ldots,q_m\in\PP^1\setminus\{p\}\), together with
frames for the corresponding subspaces, determine the determinant \(D(z)\) of
the blocks \(U_0,U_z,U_{q_3},\ldots,U_{q_m}\).  For nonzero \(z\) sufficiently
close to zero, these \(m\) fibers form a direct sum, so \(D(z)\neq0\).  This
one-dimensional slice meets \(\Delta_{12}\) transversely and meets no other diagonal.
Proposition~\ref{prop:detcriterion} therefore gives
\(\operatorname{ord}_{z=0}D(z)=e\).

We denote the wedge of the two colliding blocks by \(\Omega(z)\).  Explicitly,
$
        \Omega(z) \text{ is }
        \bigwedge_{\alpha=1}^{s}\bigwedge_{a=0}^{n_\alpha-1}u_\alpha^{(a)}(0)
        \wedge
        \bigwedge_{\alpha=1}^{s}\bigwedge_{a=0}^{n_\alpha-1}u_\alpha^{(a)}(z)
$ with \(1\leq\alpha\leq s\) and
\(0\leq a\leq n_\alpha-1\).
For every \(\alpha\) and \(0\leq a\leq n_\alpha-1\), Taylor expansion gives
\(u_\alpha^{(a)}(z)=\sum_{b\geq0}z^b u_\alpha^{(a+b)}(0)/b!\).
Expanding \(\Omega(z)\) by multilinearity, the degree in \(z\) of a typical monomial term is
\(\sum_{\alpha=1}^s\sum_{a=0}^{n_\alpha-1}b_{\alpha,a}\) for some integers \(b_{\alpha,a}\geq0\).

For a fixed \(\alpha\), the first block already contains
\(u_\alpha(0),u_\alpha'(0),\ldots,u_\alpha^{(n_\alpha-1)}(0)\).  For the
corresponding wedge term to be nonzero, the integers \(a+b_{\alpha,a}\) must be
distinct and must avoid \(0,1,\ldots,n_\alpha-1\).  Their smallest possible sum is
\(n_\alpha+(n_\alpha+1)+\cdots+(2n_\alpha-1)\).  It follows that
\(\sum_{a=0}^{n_\alpha-1}b_{\alpha,a}\) is at least the difference between
\(\sum_{j=n_\alpha}^{2n_\alpha-1}j\) and
\(\sum_{a=0}^{n_\alpha-1}a\).  This difference is \(n_\alpha^2\).
Summing over \(\alpha\), every nonzero term of \(\Omega(z)\) has order at least
\(\sum_\alpha n_\alpha^2\).  Wedging with the fixed frames of
\(U_{q_3},\ldots,U_{q_m}\) and applying a volume form cannot decrease the order of
vanishing; cancellation can only increase it.  Hence
\(\operatorname{ord}_{z=0}D(z)\geq\sum_{\alpha=1}^s n_\alpha^2\).
Since \(\operatorname{ord}_{z=0}D(z)=e\), the assertion follows.
\end{proof}

\begin{corollary}\label{cor:exponent-bounds}
Every rank-\(k\) \(m\)-unisolvent system satisfies
\(k\leq e\leq\left\lfloor\frac{\binom{mk}{k}-1}{m-1}\right\rfloor\), together with
the parity relation \(e\equiv k^2\pmod 2\).
In particular, only finitely many diagonal exponents occur for fixed \(m\) and \(k\).
\end{corollary}

\begin{proof}
Proposition~\ref{prop:lower-bound} gives
\(e\geq\sum_\alpha n_\alpha^2\geq\sum_\alpha n_\alpha=k\) and
Proposition~\ref{prop:detcriterion} gives the parity.  By
Theorem~\ref{thm:plucker-vandermonde}, the vector space underlying the Pl\"ucker
span has dimension \(e(m-1)+1\).  It is a subspace of \(\bigwedge^k V\), whose
dimension is \(\binom{mk}{k}\).  This gives the upper bound.
\end{proof}

The lower endpoint admits a complete classification.

\begin{proposition}\label{prop:minimal-exponent}
For \(m\geq2\), \(k\geq1\) and a complex vector space \(V\) of dimension
\(km\), every \(m\)-direct-sum curve \(f:\PP^1\to\Gr(k,V)\) with diagonal
exponent \(e=k\) is linearly equivalent to the standard
rank-\(k\) \(m\)-direct-sum curve. Equivalently, its associated coherent
system is isomorphic to
\(L_m^{\oplus k}=(\OO_{\PP^1}(m-1)^{\oplus k},
H^0(\PP^1,\OO_{\PP^1}(m-1)^{\oplus k}))\).
\end{proposition}

\begin{proof}
Let   $(E,W)$  be the  coherent system associated with $f$.  Suppose the bundle \(E\) splits as
\(E\simeq\bigoplus_{i=1}^k\OO_{\PP^1}(a_i)\). By Lemma \ref{lem:unisolvent-splitting-lower-bound}, we have  \(a_i\geq m-1\).  By
Proposition~\ref{prop:detcriterion}, there are identities
\(\deg E=e(m-1)=k(m-1)\). Hence \(\sum_i a_i\text{ is }k(m-1)\). So \(a_i\text{ is }m-1\) for every $i$.  Thus we have an isomorphism
\((E,V^\vee)\simeq L_m^{\oplus k}\).
\end{proof}

\begin{example} \label{par:split-hermite-models}
The rank-one Vandermonde system is
\(L_m=(\OO_{\PP^1}(m-1),H^0(\PP^1,\OO_{\PP^1}(m-1)))\), and the rank-\(k\)
split system is \(L_m^{\oplus k}\).  Its exponent is \(k\).  The rank-\(k\) Hermite system is
\(H_{m,k}=(J^{k-1}\OO_{\PP^1}(mk-1),H^0(\PP^1,\OO_{\PP^1}(mk-1)))\),
where evaluation sends a section to its \((k-1)\)-jet.  If \(D\) is reduced of
degree \(m\), the kernel of this evaluation is
\(H^0(\OO_{\PP^1}(mk-1)(-kD))=H^0(\OO_{\PP^1}(-1))=0\).  Since the source and
target have dimension \(mk\), the Hermite system is \(m\)-unisolvent.  The jet
filtration gives
\(\deg J^{k-1}\OO_{\PP^1}(mk-1)=k^2(m-1)\), so its exponent is \(k^2\).
\end{example}

We next give \(e\) a global projective interpretation by computing the degree of the
variety swept out by the moving subspaces.

\subsection{The swept variety}\label{subsec:swept-variety}

For an \(m\)-direct-sum curve \(f:\PP^1\to\Gr(k,V)\), the swept variety of \(f\) is
\(X_f=\bigcup_{t\in\PP^1}\PP(U_t)\subset\PP(V)\). If \(e\) is the diagonal exponent,
then \(\deg E=e(m-1)\).
We use the convention that \(\PP(F)\) parametrizes one-dimensional subspaces in
the fibers of a vector bundle \(F\), so
\(\OO_{\PP(F)}(-1)\subset\pi^*F\) is the tautological line subbundle.  The
inclusion \(S\subset V\otimes\OO_{\PP^1}\) induces a morphism
\(\nu:\widetilde X_f:=\PP_{\PP^1}(S)\to\PP(V)\).  It sends \((t,[u])\) to
\([u]\) and its image is \(X_f\).

\begin{proposition}\label{prop:swept-normalization-degree}
The morphism \(\nu:\widetilde X_f\to X_f\) is the normalization.  The variety
\(X_f\) is irreducible and nondegenerate of dimension \(k\) and one has
\(\deg X_f=\deg E=e(m-1)\).
\end{proposition}

\begin{proof}
If \(s\neq t\), the direct-sum condition implies
\(U_s\cap U_t=0\).  Thus \(\PP(U_s)\) and \(\PP(U_t)\) are disjoint and \(\nu\)
is injective on geometric points.  Since \(\nu\) is proper and quasi-finite, it is
finite and its degree is $1$.   Thus \(\nu\) is
birational onto its image.  The projective bundle \(\widetilde X_f\) is smooth and
therefore normal. So \(\nu\) is the normalization morphism of \(X_f\).
It follows that \(X_f\) is irreducible and
\(\dim X_f=\dim\widetilde X_f=k\).

Any \(m\) distinct subspaces \(U_{t_1},\ldots,U_{t_m}\) span \(V\), so \(X_f\)
is not contained in a hyperplane.  With \(\xi=c_1(\OO_{\PP(S)}(1))\), since
\(\nu^*\OO_{\PP(V)}(1)=\OO_{\PP(S)}(1)\) and \(\nu\) has generic degree one, the
projective-bundle formula in the convention above gives
\(\deg X_f=\int_{\PP(S)}\xi^k\).  The projective-bundle relation identifies this
integral with \(\int_{\PP^1}(-c_1(S))=\deg E\).
The last equality in the stated degree formula follows from
Proposition~\ref{prop:detcriterion}.
\end{proof}

\begin{corollary}\label{cor:swept-minimal-degree}
 The exponent \(e\text{ equals } k\) if and only if, up to a linear change of
coordinates, \(f\) is the standard split curve \(L_m^{\oplus k}\).  In this
case, \(X_f\text{ is isomorphic to }\PP^1\times\PP^{k-1}\) and
$\OO_{X_f}(1)$ is linearly equivalent to $\OO_{\PP^1\times\PP^{k-1}}(m-1,1)$.
So \(X_f\) is the balanced rational normal scroll of degree \(k(m-1)\).
\end{corollary}

\begin{proof}
  If
\(e\text{ is }k\), Proposition~\ref{prop:minimal-exponent} identifies \(f\) with the standard
split curve.  In this model,
$S$  is isomorphic to $\OO_{\PP^1}(1-m)\otimes\CC^k$ and \(\nu\) is given by the complete
linear system \(\OO_{\PP^1\times\PP^{k-1}}(m-1,1)\), which is a closed immersion.
Hence \(X_f\) is the stated Segre--Veronese scroll.  Conversely, the standard split
curve has exponent \(k\).
\end{proof}

\subsection{Collisions and Hilbert interpolation}\label{subsec:hilbert-collisions}

The definition of an \(m\)-unisolvent system tests evaluation only on reduced
divisors, which form the open locus of distinct points in the space of degree-\(m\)
effective divisors.  Passing to the Hilbert scheme is the natural extension across its
boundary: a collision is represented by the corresponding nonreduced length-\(m\)
subscheme, and the universal family retains a canonical evaluation map.

The Hilbert scheme is
\(\mathcal H_m:=\operatorname{Hilb}^m(\PP^1)=\Sym^m(\PP^1)
=|\OO_{\PP^1}(m)|\simeq\PP^m\).  The universal effective divisor is
\(\mathcal Z_m\subset\PP^1\times\mathcal H_m\), with projections \(p:\PP^1\times\mathcal H_m\to\PP^1\) and \(q:\PP^1\times\mathcal H_m\to\mathcal H_m\); we write
\(q_{\mathcal Z}:=q|_{\mathcal Z_m}\).  Under the identification of \(\mathcal H_m\)
with the complete linear system, the ideal of \(\mathcal Z_m\) is
\(\mathcal I_{\mathcal Z_m}\simeq
p^*\OO_{\PP^1}(-m)\otimes q^*\OO_{\mathcal H_m}(-1)\).  For a vector bundle \(E\)
of rank \(k\) on \(\PP^1\), its length-\(m\) tautological bundle is
\(E^{[m]}:=q_{\mathcal Z*}(p^*E|_{\mathcal Z_m})\).
The morphism \(q_{\mathcal Z}\) is finite and flat of degree \(m\).  Thus
\(E^{[m]}\) is locally free of rank \(mk\), its formation commutes with arbitrary
base change and its fiber at an effective divisor \(D\), not necessarily reduced, is
canonically isomorphic to \(H^0(E|_D)\).  A coherent system \((E,W)\) with \(\dim W=mk\) therefore
has the universal Hilbert evaluation morphism
\(\epsilon_{E,W}:W\otimes\OO_{\mathcal H_m}\to E^{[m]}\).

The following resolution is the vector-bundle form of the classical
Schwarzenberger construction; see \cite{Schwarzenberger64,RaicuSam22} for the
line-bundle construction and its modern representation-theoretic treatment.

\begin{proposition}\label{prop:tautological-resolution}
For an \(m\)-unisolvent system \((E,W)\) of rank \(k\), degree \(d\) and
diagonal exponent \(e\), the integer
\(b:=h^0(E(-m))=d-k(m-1)=(e-k)(m-1)\)
is nonnegative.  Moreover, \(E^{[m]}\) has the natural two-term resolution
\(0\longrightarrow
H^0(E(-m))\otimes\OO_{\mathcal H_m}(-1)\xrightarrow{\ \mu\ }
H^0(E)\otimes\OO_{\mathcal H_m}\longrightarrow E^{[m]}\longrightarrow0\).
Its determinant line bundle satisfies
\(\det E^{[m]}\simeq
\det H^0(E)\otimes\det H^0(E(-m))^{-1}\otimes\OO_{\mathcal H_m}(b)\).
In particular, \(c_1(E^{[m]})\text{ is }bH\), where \(H\text{ is }c_1(\OO_{\mathcal H_m}(1))\).
\end{proposition}

\begin{proof}
Tensoring the universal divisor sequence by \(p^*E\) gives an exact sequence  \(0\to p^*E(-m)\otimes q^*\OO_{\mathcal H_m}(-1)\to
p^*E\to p^*E|_{\mathcal Z_m}\to 0\). (*)
  By the projection formula, we have $R^kq_*(p^*E(-m)\otimes q^*\OO_{\mathcal H_m}(-1))=H^k(E(-m))\otimes \OO_{\mathcal H_m}(-1)$ and $R^kq_*(p^*E)=H^k(E)\otimes \OO_{\mathcal H_m}$.  Lemma~\ref{lem:unisolvent-splitting-lower-bound} gives
\(H^1(E)=H^1(E(-m))=0\). Applying \(R^kq_*\) to the exact sequence $(*)$, we get
 the stated resolution.  The formula for \(b\) follows from
Riemann--Roch on \(\PP^1\). Taking determinants gives the stated determinant-line
isomorphism.
\end{proof}

The arrow \(\mu\) in this resolution is multiplication by the
universal binary form.  At \([F]\in|\OO_{\PP^1}(m)|\), it specializes to
\(H^0(E(-m))\otimes\langle F\rangle\xrightarrow{\,\cdot F\,}H^0(E)\) where $\langle F\rangle$ is the fiber of $\OO_{\mathcal H_m}(-1)$ at $[F]$.
The discriminant hypersurface \(\Delta_m\subset\mathcal H_m\simeq\PP^m\) parametrizes
nonreduced effective divisors, or equivalently binary forms with a
multiple root.  It is irreducible and has degree \(2(m-1)\).

\begin{theorem}\label{thm:hilbert-discriminant}
For an \(m\)-unisolvent system \((E,W)\) of rank \(k\) and diagonal exponent \(e\),
the determinant of the universal Hilbert evaluation satisfies the scheme-theoretic
identity \(\divisor(\det\epsilon_{E,W})=\frac{e-k}{2}\,\Delta_m\).

\end{theorem}

\begin{proof}
Over \(\mathcal H_m\setminus\Delta_m\), a point represents a reduced divisor of
degree \(m\), so unisolvence says precisely that \(\epsilon_{E,W}\) is an
isomorphism.  Hence \(\divisor(\det\epsilon_{E,W})\) is supported on the irreducible
hypersurface \(\Delta_m\), so
\(\divisor(\det\epsilon_{E,W})=a\Delta_m\) for some nonnegative integer \(a\).  By
Proposition~\ref{prop:tautological-resolution}, the determinant is a nonzero section
of a line bundle of numerical class \(bH=(e-k)(m-1)H\).  Comparing degrees with
\(\deg\Delta_m=2(m-1)\) gives \(2a(m-1)=(e-k)(m-1)\), and hence
\(a=(e-k)/2\).

\end{proof}
\begin{corollary}
\label{cor:hilbert-isomorphism-locus}
Let \((E,W)\) be an \(m\)-unisolvent system of rank \(k\) and
diagonal exponent \(e\). The isomorphism locus of
$
\epsilon_{E,W}\colon
W\otimes\mathcal O_{H_m}\longrightarrow E^{[m]}
$
is
$H_m
\text{ if }e=k$
and is exactly
$H_m\setminus\Delta_m
\text{ if }e>k.$
In particular, \(e=k\) if and only if the restriction map
$W\longrightarrow H^0(E|_D)$
is an isomorphism for some nonreduced closed subscheme
\(D\subset\mathbb P^1\) of length \(m\).
\end{corollary}

\begin{proof}
By \(m\)-unisolvence, the morphism \(\epsilon_{E,W}\) is an
isomorphism over \(H_m\setminus\Delta_m\). By Theorem~\ref{thm:hilbert-discriminant}, we have
$
\operatorname{div}(\det\epsilon_{E,W})
=
\frac{e-k}{2}\Delta_m.
$
If \(e=k\), the determinant is nowhere vanishing, whereas if \(e>k\),
it vanishes at every point of \(\Delta_m\). The assertions follow.
\end{proof}

 \section{Rank-two classification and geometric applications}\label{sec:rank-two-applications}
\subsection{Rank-two classification}
The exponent theory of the preceding section identifies the extremal models, but a complete
classification also requires control of the local Perkinson partition. Rank two is the first genuinely
Grassmannian case and has only the two partitions \((2)\) and \((1,1)\). This finite alternative makes a
complete classification possible and turns the numerical theory into concrete projective geometry.
The rank-two classification is therefore the technical core of the paper. Theorem~\ref{thm:single-string}
treats the partition \((2)\), identifying the curve with the Hermite model \(H_{m,2}\); only the
partition \((1,1)\) remains.

The geometry behind the proof is already visible in \(\Gr(2,4)\). For a type-\((1,1)\) curve with
local frame \(U_t=\Span\{u_1(t),u_2(t)\}\), one has
\(U_t^{(1)}=\Span\{u_1(t),u_2(t),u_1'(t),u_2'(t)\}=V\) at a general point.
Writing \(\omega(t):=u_1(t)\wedge u_2(t)\), Taylor expansion at a collision gives
\(\omega(t)\wedge\omega(t+\varepsilon)=\varepsilon^2u_1(t)\wedge u_2(t)\wedge
u_1'(t)\wedge u_2'(t)+O(\varepsilon^3)\).
Thus the diagonal order is \(2\) and Proposition~\ref{prop:minimal-exponent} proves the
type-\((1,1)\) case in \(\Gr(2,4)\).

The first instructive higher-dimensional case is \(\Gr(2,6)\). For \(V=\CC^6\) and
\(y\in\PP^1\), the rational projection \(\pi_y:\Gr(2,V)\dashrightarrow\Gr(2,V/U_y)\) has center
\(U_y\). For distinct
\(x_1,x_2\in\PP^1\setminus\{y\}\), the \(3\)-direct-sum property gives
\(U_y\oplus U_{x_1}\oplus U_{x_2}=V\).  After quotienting by \(U_y\), one obtains
\(\overline U_{x_1}\oplus\overline U_{x_2}=V/U_y\).  Thus the projected curve has the
\(2\)-direct-sum property away from \(y\).

The projected curve extends across \(y\), but its value there is not \(U_y/U_y\). In a coordinate
with \(y=0\), a frame \(U_t=\Span\{u_1(t),u_2(t)\}\) satisfies
\(u_i(t)=u_i(0)+tu_i'(0)+O(t^2)\), and hence
\(\overline{u_i(t)}=t\,\overline{u_i'(0)}+O(t^2)\) modulo \(U_0\). Its value at \(y\) is therefore
\(U_y^{(1)}/U_y\), where \(U_y^{(1)}\) is spanned by
\(\{u_1(0),u_2(0),u_1'(0),u_2'(0)\}\). An inductive projection
argument consequently requires the transversality \(U_y^{(1)}\oplus U_q=V\).  This is the
geometric reason for Lemma~\ref{lem:rank-two-osculating}.

\begin{theorem}\label{thm:rank-two-direct-sum}
Every \(m\)-direct-sum curve \(f:\PP^1\to\Gr(2,V)\) with \(\dim V=2m\) whose normal form has type
\((1,1)\) at a general point admits a decomposition \(V=V_1\oplus V_2\), where \(\dim V_i=m\).
Moreover, there are rational normal curves \(C_i\subset\PP(V_i)\simeq\PP^{m-1}\) of degree \(m-1\)
such that \(f(t)=\ell_1(t)\oplus\ell_2(t)\), where \(\ell_i(t)\subset V_i\) is the line represented by
\(C_i(t)\).
\end{theorem}

\begin{proof}
The pullback of the tautological subbundle is denoted by
\(S\subset V\otimes\OO_{\PP^1}\), with \(U_t=f(t)=S_t\subset V\). For a local holomorphic frame
\((a(t),b(t))\), the first osculating subspace of \(U_t\) is
\(U_t^{(1)}:=\Span\{a(t),b(t),a'(t),b'(t)\}\subset V\). This definition is independent of the frame.
Since the normal form has type \((1,1)\) on a
nonempty open subset, a general point \(p\) satisfies \(\dim U_p^{(1)}=4\).

\begin{lemma}\label{lem:rank-two-osculating}
For a general point \(p\in\PP^1\), there are general points \(q_3,\ldots,q_m\in\PP^1\) such that
\(U_p^{(1)}\oplus U_{q_3}\oplus\cdots\oplus U_{q_m}=V\).
\end{lemma}

\begin{proof}[Proof of the lemma]
For \(m=2\), the assertion is immediate because \(\dim V=4=\dim U_p^{(1)}\). The case \(m\geq3\)
remains.
 The following elementary fact will be used.

\begin{fact}\label{fact:two-dimensional}
Let $B$ be a two-dimensional vector space and let $\eta(t)\in B$ be holomorphic on a small disc $\Delta$, with
$\eta(0)\neq0$.  Suppose that for all sufficiently small distinct nonzero $t,s$, the vectors
$\eta(t)$ and $\eta(s)$ are linearly independent.  Then by fixing a volume element of $\wedge^2 B$,  $\eta(t)\wedge\eta(s)$ can be viewed as $(t-s)v(t,s)$, where
$v(t,s)$ is a unit.
\end{fact}

We now prove the fact. By shrinking the disc and without loss of generality,  we can write $\eta(t)=c(t)(1,z(t))$, where $c(t)$ is a unit.  Then we have $z(t)\neq z(s)$ for nonzero $t\neq s$.  If $z(t)=z(0)+\alpha t^r+O(t^{r+1})$ with $r\geq2$, then $z(t)$ is $r$-ramified on any small disc around $0$, which is a
contradiction.  Hence $r=1$, so $(z(s)-z(t))/(s-t)$ is a unit.

We first handle the case $m=3$.  Let $p=0$ and $A=U_0^{(1)}$, so $\dim A=4$.  We claim that $A\cap U_t=0$ for
$t\neq0$ near $0$. Suppose otherwise. The quotient map
$V \rightarrow V/A
$
induces a holomorphic bundle morphism
$\Phi:S\rightarrow (V/A)\otimes\mathcal{O}_{\Delta}.
$
Since \(\det \Phi\) vanishes at points accumulating at \(0\), it vanishes identically. Hence the generic rank of \(\Phi\) is at most one. It cannot be zero, since otherwise \(U_t\subset A\) for general \(t\), contradicting the \(3\)-direct-sum condition. After shrinking \(\Delta\), the kernel \(\ker \Phi\) has rank one on
$\Delta^{*}:=\Delta\setminus\{0\}.
$
Its saturated extension in \(\mathcal{S}\) is a line subbundle on \(\Delta\). Choose a nowhere-vanishing holomorphic section \(x(t)\) of this line subbundle over $\Delta$. Then
$x(t)\in A\cap U_t
$
for every \(t\neq 0\). Complete \(x(t)\) to a holomorphic frame \(x(t),y(t)\) of \(\mathcal{S}\). Let $e_1=x(0)$ and $e_2=y(0)$ and choose a
decomposition $V=A\oplus B$ with $\dim B=2$. Let $\bar y(t)$ be the image of $y(t)$ in $B$.

For a local frame $a(t),b(t)$ of $U_t$, the determinant
$
        \Delta(t_1,t_2,t_3)=a(t_1)\wedge b(t_1)\wedge a(t_2)\wedge b(t_2)
        \wedge a(t_3)\wedge b(t_3)
$
 can be written as
$
        \Delta(t_1,t_2,t_3)=(t_1-t_2)^e(t_1-t_3)^e(t_2-t_3)^e h(t_1,t_2,t_3)
$
for a unit $h$.  Replacing the frame by $(e_1,e_2)$ at $0$ and by $(x(t),y(t))$ at $t$ changes the
determinant only by units.  Thus
$
        D(t,s):=e_1\wedge e_2\wedge x(t)\wedge y(t)\wedge x(s)\wedge y(s)
        =t^es^e(t-s)^e h_0(t,s)
$
for a unit $h_0$.

As $V=A\oplus B$, the same determinant factors as
$
        D(t,s)=\pm\bigl(e_1\wedge e_2\wedge x(t)\wedge x(s)\bigr)
        \bigl(\bar y(t)\wedge\bar y(s)\bigr).
$
Let $\bar x(t)$ be the image of $x(t)$ in $A/U_0$.  Since $x(0)\in U_0$, write $\bar x(t)=t^a\xi(t)$ with
$\xi(0)\neq0$.  The integer $a$ is finite, because $U_0\cap U_t=0$ for $t\neq0$ by the 3-direct sum
condition.  If $a\geq2$, then $x'(0)\in U_0$,  which contradicts the fact $A=U_0^{(1)}=\Span\{x(0),y(0),x'(0),y'(0)\}$.  Hence $a$ is $1$.  Applying Fact~\ref{fact:two-dimensional} to
$A/U_0$, we have
$
        e_1\wedge e_2\wedge x(t)\wedge x(s)=ts(t-s)u(t,s)
$
for a unit $u$.

Since $y(0),y'(0)\in A$, the image of $y(t)$ in $B=V/A$ has order at least two: $\bar y(t)=t^\rho\eta(t)$
with $\rho\geq2$ and $\eta(0)\neq0$.  For nonzero $t,s$ with $t\neq s$, the vectors $\bar y(t)$ and $\bar
y(s)$ are linearly independent; otherwise one would obtain a nontrivial relation among the direct
summands $U_0$, $U_t$ and $U_s$.  By Fact~\ref{fact:two-dimensional}, we have
$
        \bar y(t)\wedge\bar y(s)=t^\rho s^\rho(t-s)v(t,s)
$
for a unit $v$.  Consequently, there is an identity
$
        D(t,s)=t^{\rho+1}s^{\rho+1}(t-s)^2w(t,s)
$
with $w$ a unit.  Comparing with $t^es^e(t-s)^e h_0(t,s)$,  $e$ must be $2$ and must also equal $\rho+1$. Note
$\rho+1\geq3$,  we get a contradiction.  So $A\cap U_t=0$ for
$t\neq0$ near $0$, and we prove the $m=3$ case.
Now assume that \(m>3\). Fix a general point \(p\) and put \(A:=U_p^{(1)}\), so that \(\dim A=4\).

We first choose general points \(q_4,\ldots,q_m\) such that, with \(B':=U_{q_4}\oplus\cdots\oplus U_{q_m}\), one has \(A\cap B'=0\). We construct these points inductively. Suppose that, for some \(0\le r\le m-4\), we have already chosen \(r\) pairwise distinct general points \(q_4,\ldots,q_{r+3}\), all different from \(p\), such that \(B_r:=U_{q_4}\oplus\cdots\oplus U_{q_{r+3}}\) satisfies \(A\cap B_r=0\). Here \(B_0=0\). Since the original curve is \(m\)-direct, every collection of at most \(m\) distinct fibers is a direct sum; in particular, \(\dim B_r=2r\). Hence \(\dim(A\oplus B_r)=4+2r\).

We claim that a general further fiber \(U_q\) is disjoint from \(A\oplus B_r\). Suppose otherwise. Then, for general \(q\), one has \(U_q\cap(A\oplus B_r)\neq0\). Choose \(m-r-1\) further general pairwise distinct points \(t_1,\ldots,t_{m-r-1}\), different from \(p,q_4,\ldots,q_{r+3}\). By the preceding assumption, each \(U_{t_i}\) meets the fixed space \(A\oplus B_r\) nontrivially. Consequently, after adjoining these fibers one at a time, each of them contributes at most one new dimension. Therefore
\(\dim(A+B_r+U_{t_1}+\cdots+U_{t_{m-r-1}})\le 4+2r+(m-r-1)=m+r+3\).
Since \(r\le m-4\), we have \(m+r+3\le 2m-1\).

On the other hand, the space \(A+B_r+U_{t_1}+\cdots+U_{t_{m-r-1}}\) contains the sum of the \(m\) distinct fibers \(U_p,U_{q_4},\ldots,U_{q_{r+3}},U_{t_1},\ldots,U_{t_{m-r-1}}\). By the \(m\)-direct-sum condition, these \(m\) fibers form a direct sum equal to \(V\), and hence their span has dimension \(2m\). This contradicts the preceding upper bound \(2m-1\). Thus a general further fiber \(U_q\) satisfies \(U_q\cap(A\oplus B_r)=0\). Taking such a \(q\) and setting \(B_{r+1}:=B_r\oplus U_q\), we obtain \(A\cap B_{r+1}=0\). Proceeding inductively, we obtain general points \(q_4,\ldots,q_m\) such that \(B'=U_{q_4}\oplus\cdots\oplus U_{q_m}\) and \(A\cap B'=0\).
Now set \(\overline V:=V/B'\), and let \(\pi:V\to\overline V\) be the quotient map. Since \(\dim B'=2(m-3)\), we have \(\dim\overline V=2m-2(m-3)=6\). For \(t\notin\{q_4,\ldots,q_m\}\), let \(\overline U_t:=\pi(U_t)\). The restriction \(U_t\to\overline U_t\) is injective.  Thus we have \(\dim\overline U_t=2\).

Moreover, if \(t_1,t_2,t_3\) are pairwise distinct points lying in   \(\mathbb{P}^1\setminus \{q_4,\ldots,q_m\}\), then the \(m\)-direct-sum condition gives \(V=U_{t_1}\oplus U_{t_2}\oplus U_{t_3}\oplus B'\). Passing to the quotient by \(B'\), we obtain \(\overline V=\overline U_{t_1}\oplus\overline U_{t_2}\oplus\overline U_{t_3}\). Hence, near \(p\), the projected family \(t\mapsto\overline U_t\) is a \(3\)-direct-sum family in \(\operatorname{Gr}(2,\overline V)\).
Since \(A\cap B'=0\), the restriction of \(\pi\) to \(A\) is injective. If \(a(t),b(t)\) is a local holomorphic frame of \(U_t\) near \(p\), then \(\overline U_p^{(1)}=\operatorname{Span}\{\pi(a(p)),\pi(b(p)),\pi(a'(p)),\pi(b'(p))\}=\pi(A)\). Thus \(\dim\overline U_p^{(1)}=4\), and \(\pi\) identifies \(A\) with \(\overline U_p^{(1)}\).
The proof of the case \(m=3\) is local near the chosen point and uses only the \(3\)-direct-sum property. We may therefore apply it to the projected family. It follows that, for a general point \(q_3\), one has \(\overline U_p^{(1)}\oplus\overline U_{q_3}=\overline V\). Equivalently, \(\pi(A)+\pi(U_{q_3})=\overline V\), and hence \(A+U_{q_3}+B'=V\). Finally, \(\dim A+\dim U_{q_3}+\dim B'=4+2+2(m-3)=2m=\dim V\). Therefore the above sum is direct, and we obtain \(U_p^{(1)}\oplus U_{q_3}\oplus\cdots\oplus U_{q_m}=V\). This proves the lemma for \(m>3\).
\end{proof}

We now compute the diagonal exponent.  Choose $p,q_3,\ldots,q_m$ as in
Lemma~\ref{lem:rank-two-osculating} and choose a local frame $u_1(t),u_2(t)$ of the tautological
subbundle near $p$.  Let $t_2=p+\varepsilon$.  Taylor expansion gives
$u_1(p)\wedge u_2(p)\wedge u_1(p+\varepsilon)\wedge u_2(p+\varepsilon)
  \wedge u_1(q_3)\wedge u_2(q_3)\wedge\cdots\wedge u_1(q_m)\wedge u_2(q_m)
=\varepsilon^2 u_1(p)\wedge u_2(p)\wedge u_1'(p)\wedge u_2'(p)
  \wedge u_1(q_3)\wedge u_2(q_3)\wedge\cdots\wedge u_1(q_m)\wedge u_2(q_m)
  +O(\varepsilon^3).$
The coefficient of $\varepsilon^2$ is nonzero by Lemma~\ref{lem:rank-two-osculating}.  Hence the
vanishing order along a general point of any diagonal is exactly $2$, so the diagonal exponent is $e=2$.
Applying Proposition~\ref{prop:minimal-exponent} with $k=2$, the associated coherent system is $L_m\oplus
L_m$.
\end{proof}

\begin{proposition}\label{prop:low-rank}
For every \(m\geq2\), the following statements hold.
\begin{enumerate}[label=\textnormal{(\arabic*)},nosep]
\item Up to isomorphism, the unique rank-one \(m\)-unisolvent system is \(L_m\), whose bundle is
\(\OO_{\PP^1}(m-1)\) and whose space of sections is
\(H^0\bigl(\PP^1,\OO_{\PP^1}(m-1)\bigr)\).
\item Up to isomorphism, the rank-two \(m\)-unisolvent systems are \(L_m\oplus L_m\) and \(H_{m,2}\),
where the bundle of \(H_{m,2}\) is \(J^1\OO_{\PP^1}(2m-1)\) and its space of sections is
\(H^0\bigl(\PP^1,\OO_{\PP^1}(2m-1)\bigr)\).
\end{enumerate}
\end{proposition}

\begin{proof}
The rank-one case is the rational normal curve in \(\PP^{m-1}\). In rank two, the normal-form partition
is either \((2)\) or \((1,1)\). The assertion follows from Theorems~\ref{thm:single-string}
and~\ref{thm:rank-two-direct-sum}.
\end{proof}

\subsection{Geometric applications}\label{subsec:geometric-applications}

The rank-two classification and the preceding general theory have three
consequences in projective geometry.
\paragraph{Critically spanning line rulings.}
\begin{theorem}\label{thm:critical-spanning-scroll}
 Let $m\geq2$, let $V$ be a complex vector space of dimension $2m$ and let
$X$ be a positive-dimensional  projective subvariety of $\Gr(2,V)$. Assume that $Y\subset \PP(V)$ is a closed subvariety  ruled  by lines over $X$. Suppose  $\langle\ell_{x_1},\ldots,\ell_{x_m}\rangle=\PP(V)$ for every
$m$ pairwise distinct closed points $x_1,\ldots,x_m\in X$ where $\ell_{x_i}$ denotes the line over $x_i$.
Then
$X$ is a smooth rational curve  and $Y$ is an irreducible nondegenerate surface.  Moreover, up to projective
equivalence, exactly one of the following cases occurs.
\begin{enumerate}[label=\textup{(\arabic*)},nosep]
\item  $Y$ is the balanced rational normal scroll $S_{m-1,m-1}$.  The isomorphism
$Y\simeq\PP^1\times\PP^1$ identifies $\OO_Y(1)$ with
$\OO_{\PP^1\times\PP^1}(m-1,1)$.  In particular, $Y$ is smooth
and $\deg Y\text{ is }2m-2$.
\item There is a degree-$(2m-1)$ rational normal embedding
$c:\PP^1\to\PP(V)$, with image $R$, such that
$\ell_t=T_{c(t)}R$.  In this case
$Y=\operatorname{Tan}(R)$ and its normalization is
$q:\PP^1\times\PP^1\to Y$, with
$q^*\OO_Y(1)\simeq\OO_{\PP^1\times\PP^1}(2m-2,1)$.  Moreover,
$\deg Y\text{ is }4m-4$, the variety $Y$ is nonnormal and
$\operatorname{Sing}(Y)_{\mathrm{red}}=R$.
\end{enumerate}
Both families satisfy the stated spanning condition.  Consequently, $Y$ is
smooth if and only if it is normal, and these conditions hold precisely in
the first case.
\end{theorem}

\begin{proof}
The ruling defines a morphism \(f:X\to\Gr(2,V)\), where
\(\PP(U_x)=\ell_x\). Since the dimensions of any \(m\) fibers add up to
\(\dim V\), the spanning hypothesis is equivalent to the
\(m\)-direct-sum property. Theorem~\ref{source} gives \(X\simeq\PP^1\),
and Proposition~\ref{prop:low-rank} identifies the associated system with
\(L_m\oplus L_m\) or \(H_{m,2}\). Put \(S=f^*\mathcal U\). By
Proposition~\ref{prop:swept-normalization-degree}, the natural morphism
\(q:\PP_{\PP^1}(S)\to Y\) is the normalization, and \(Y\) is an irreducible
nondegenerate surface with \(\deg Y=e(m-1)\).

In the split case, \(S\simeq\OO_{\PP^1}(1-m)\otimes\CC^2\), so
\(\PP_{\PP^1}(S)\simeq\PP^1\times\PP^1\). With the convention of
Section~\ref{subsec:swept-variety}, its tautological line subbundle is
\(\OO_{\PP^1\times\PP^1}(1-m,-1)\). Therefore
\(q^*\OO_Y(1)\simeq\OO_{\PP^1\times\PP^1}(m-1,1)\).
The morphism \(q\) is defined by this complete linear system, which is
very ample. Thus it is an isomorphism onto the balanced rational normal
scroll \(S_{m-1,m-1}\), and \(Y\) is smooth of degree \(2m-2\), as also
follows from Corollary~\ref{cor:swept-minimal-degree}.

In the Hermite case, set \(d:=2m-1\) and choose a two-dimensional vector
space \(B\). After a linear change of coordinates, we may write
\(V=\Sym^d B\), with the rational normal curve parametrized by
\(c([l])=[l^d]\). Its affine tangent space at \([l^d]\) is
\(U_{[l]}=l^{d-1}B\). Multiplication by \(l^{d-1}\) identifies
\(S\) with \(\OO_{\PP(B)}(1-d)\otimes B\). Consequently,
\[
q:\PP(B)\times\PP(B)\longrightarrow Y,
\qquad ([l],[w])\longmapsto[l^{d-1}w],
\qquad
q^*\OO_Y(1)\simeq\OO_{\PP^1\times\PP^1}(d-1,1).
\]
The image is \(\operatorname{Tan}(R)\), where
\(R=\{[l^d]:[l]\in\PP(B)\}\), and the polarization is
\(\OO_{\PP^1\times\PP^1}(2m-2,1)\). Its self-intersection is
\(2(d-1)=4m-4\), which gives the asserted degree because \(q\) is
birational.

It remains to determine the singular locus in the Hermite case. Since
\(d\geq3\), unique factorization of binary forms shows that
\([l^{d-1}w]\) determines both \([l]\) and \([w]\). The same argument
holds over every algebraically closed extension of \(\CC\), so \(q\) is
universally injective. Moreover, \(q^{-1}(R)\) is the diagonal. At a point
off the diagonal, choose a basis \(x,y\) of \(B\) with \(l=x\) and
\(w=y\). The two tangent directions of \(q\), modulo the line
\(\CC x^{d-1}y\), are represented by
\((d-1)x^{d-2}y^2\) and \(x^d\), and hence are independent. Thus the
restriction \(q^{-1}(Y\setminus R)\to Y\setminus R\) is unramified.
It is also finite and universally injective, so it is a closed immersion
by \cite[Tag~04XV]{StacksProject}. Since its target is reduced and it is
surjective, it is an isomorphism. In particular, \(Y\) is smooth away
from \(R\).

At a point of the diagonal, choose \(l=w=x\). The two tangent directions
are now represented by \((d-1)x^{d-1}y\) and \(x^{d-1}y\), so \(dq\) has
rank one. If \(Y\) were normal at the corresponding point of \(R\), its
normalization would be an isomorphism there, contradicting this rank
drop. Thus every point of \(R\) is nonnormal and hence singular. We
obtain \(\operatorname{Sing}(Y)_{\mathrm{red}}=R\), and \(Y\) is
nonnormal.

Finally, both the split and Hermite systems are \(m\)-unisolvent by
Example~\ref{par:split-hermite-models}. Proposition~\ref{prop:equivalence}
therefore gives the spanning condition for both families. The assertions
about smoothness and normality follow from the two cases above.
\end{proof}

For $m=2$, the spanning condition says exactly that distinct ruling lines
in $\PP^3$ are disjoint.  The two surfaces in
Theorem~\ref{thm:critical-spanning-scroll} are a smooth quadric and the
quartic tangent developable of the twisted cubic.  If the ruling is
presented by two disjoint embedded directrices, these define disjoint
sections of $\PP^1\times\PP^1\to\PP^1$.  For a polarization $\OO(a,1)$,
the graph of a degree-$d$ morphism $\PP^1\to\PP^1$ has class $dA+B$, where
$A$ and $B$ are the pullbacks of the point classes on the two factors; its
hyperplane degree is $a+d$.
Two such sections of degrees $d_1$ and $d_2$ have intersection number
$d_1+d_2$.  Disjointness therefore gives $d_1=d_2=0$, so both directrices
have degree $a$.  Here $a=1$ in the split
case and $a=2$ in the Hermite case.  Thus the possible bidegrees are
$(1,1)$ and $(2,2)$, and no genericity hypothesis is required.  This
recovers Day's two families and excludes all further bidegrees
\cite[Lemma~3.11, Theorem~3.13, and Example~3.14]{Day2025}.

\paragraph{Higher-order critical Terracini rigidity.}
Let
$X\subset\PP(V)$ be a  smooth, projective and nondegenerate
variety of dimension $n$ $  (n\ge 1)$. We set $\rho_s(n):=\binom{n+s}{n}$ with  $s\geq1$. Let
$\widehat{\operatorname{Osc}}^s_xX\subset V$ denote the affine $s$-th
osculating space at $x$.  Assume that  $\dim V\text{ is }m\rho_s(n)$ for some integer $m\ge 2$. The reduced critical
$s$-th Terracini locus is
$\mathcal T_m^{(s),\mathrm{red}}(X):=\{(x_1,\ldots,x_m)\in
\operatorname{Conf}_m(X):\widehat{\operatorname{Osc}}^s_{x_1}X+
\cdots+\widehat{\operatorname{Osc}}^s_{x_m}X\neq V\}$, where
$\operatorname{Conf}_m(X)$ is the space of ordered $m$-tuples of pairwise
distinct points.

\begin{theorem}\label{thm:higher-critical-terracini}
Suppose that $\dim V\text{ is }m\rho_s(n)$.  Then
$\mathcal T_m^{(s),\mathrm{red}}(X)$ is empty if and only if $n\text{ is }1$ and
$X\subset\PP(V)$ is projectively equivalent to the rational normal curve
$\nu_{m(s+1)-1}(\PP^1)\subset\PP^{m(s+1)-1}$.
\end{theorem}
\begin{proof}
Set $L:=\OO_X(1)$.  The rank of the target of the principal-parts morphism
$j^s:V^\vee\otimes\OO_X\to J_X^sL$ is $\rho_s(n)$ and
the image of $(j_x^s)^\vee$ is
$\widehat{\operatorname{Osc}}^s_xX$.  If the critical Terracini locus is
empty, the equality $\dim V=m\rho_s(n)$ forces $j^s$ to be surjective at
every point and the $m$ osculating spaces attached to any reduced
$m$-tuple to form a direct sum.  Hence the dual subbundle
$(J_X^sL)^\vee\subset V\otimes\OO_X$ defines a nonconstant
$m$-direct-sum morphism
$\gamma_s:X\to\Gr(\rho_s(n),V)$. By Propositions~\ref{prop:dimension-bound}
and~\ref{prop:curve-is-p1}, we have $n=1$ and $X\simeq\PP^1$.
Put $r:=s+1$.  If $v(t)$ is a local lift of the embedding, then
$\gamma_s(t)$ is spanned by $v(t),v'(t),\ldots,v^{(s)}(t)$.  So
its drill is one.  By Theorem~\ref{thm:single-string},  $X\subset\PP(V)$ is projectively equivalent to the rational normal curve
$\nu_{m(s+1)-1}(\PP^1)\subset\PP^{m(s+1)-1}$.

Conversely, set $N:=m(s+1)$ and consider the rational normal curve defined
by $|\OO_{\PP^1}(N-1)|$.  For pairwise distinct $x_1,\ldots,x_m$, the
annihilator of the sum of the affine $s$-th osculating spaces is
$H^0\!\left(\OO_{\PP^1}(N-1)\bigl(-(s+1)(x_1+\cdots+x_m)\bigr)\right)
=H^0(\OO_{\PP^1}(-1))=0$.
The osculating spaces therefore span $V$, proving the converse.
\end{proof}

For $s=1$, the theorem specializes to critical tangent-space rigidity; its
curve case is the result of Ballico and Chiantini
\cite[Proposition~7]{BallicoChiantini2024}.

\paragraph{Co-critical conormal rigidity.}

Let $X\subset\PP^N=\PP(V)$ be a  smooth, projective and
nondegenerate subvariety of dimension $n\ge 1$. We set $c:=N-n$ and assume $N+1\text{ is }mc$ ($m\ge2$).

\begin{theorem}\label{thm:co-critical-rigidity}
Suppose that $n\text{ is at least }1$ and the intersection satisfies
$T_{x_1}X\cap\cdots\cap T_{x_m}X=\varnothing$ for every $m$ pairwise
distinct points of $X$.  Then exactly one of the following occurs:
\begin{enumerate}[label=\textup{(\arabic*)},nosep]
\item $m$ is $2$, $N$ is $3$ and $X$ is the twisted cubic;
\item $m$ is $3$, $N$ is $2$ and $X$ is a smooth conic.
\end{enumerate}
Conversely, both curves satisfy the stated intersection condition.
\end{theorem}

\begin{proof}
The hypothesis is equivalent to
$\bigcap_i\widehat T_{x_i}X=0$.  Taking annihilators gives
$\sum_i(\widehat T_{x_i}X)^\perp=V^\vee$, and the equality
$mc=N+1$ makes this sum direct.  Thus $\gamma_X^\perp$ has the
$m$-direct-sum property.  Propositions~\ref{prop:dimension-bound}
and~\ref{prop:curve-is-p1} yield $n=1$ and $X\simeq\PP^1$.  Since then
$c=N-1$, the critical equality becomes $(m-1)(N-1)=2$, whose only solutions
are $(m,N)=(2,3)$ and $(m,N)=(3,2)$.

For $(m,N)=(2,3)$, the tangent lines at any two distinct points are
disjoint, so Theorem~D with $s=1$
identifies $X$ with the twisted cubic.  For $(m,N)=(3,2)$, the image of
$\gamma_X^\perp$ is an integral nondegenerate curve in the dual plane, and
the direct-sum condition says that no three of its distinct points are
collinear.  Considering the intersection number with a general line, its degree is at most two;
nondegeneracy implies that it is a smooth conic.  Projective biduality then shows that
$X$ is a smooth conic \cite{TevelevDuality}.  Conversely, distinct tangent
lines to the twisted cubic are skew, while three distinct tangent lines to
a smooth conic are not concurrent.
\end{proof}

 \section{The unisolvent category and the rank-three threshold}\label{sec:unisolvent-category}

The rank-two classification determines all simple objects of ranks at most
two. We now use these low-rank objects to study rank three. The
\(m\)-unisolvent systems form an abelian category of finite length, and
the diagonal exponent is additive in short exact sequences. These facts
restrict the exponent of a nonsimple rank-three system to \(3\) or \(5\).
For \(m=2\), both values are controlled by the rank-one and rank-two
models, whereas a simple rank-three system of exponent \(7\) exists.
Thus \(7\) is the first exponent at which a new simple object occurs in
this case.

\subsection{The abelian length category}

\begin{definition}
Fix \(m\geq2\). Let \(\mathcal U_m\) be the category whose objects are
\(m\)-unisolvent systems on \(\PP^1\). A morphism from \((E,W)\) to
\((F,Z)\) is a pair \((\varphi,\psi)\), where \(\varphi:E\to F\) is a
bundle morphism, \(\psi:W\to Z\) is a linear map and the restriction of
\(H^0(\varphi)\) to \(W\) equals \(\psi\). The zero system \((0,0)\) is
also an object.
\end{definition}

The category of coherent systems embeds in the abelian category
\(\mathcal C\) whose objects consist of a finite-dimensional vector space
\(W\), a coherent sheaf \(E\) and a structure morphism
\(W\otimes\OO_{\PP^1}\to E\); see
\cite[Section~2.5]{KingNewstead95}.

\begin{lemma}\label{lem:full-subcategory}
Inside \(\mathcal C\), the category \(\mathcal U_m\) is precisely the full
subcategory \(\mathcal D_m\) of objects whose structure morphism induces
an isomorphism \(W\to H^0(E|_D)\) for every reduced divisor \(D\) of
degree \(m\).
\end{lemma}

\begin{proof}
It is enough to show that for \((E,W)\in\mathcal D_m\), the sheaf \(E\)
is locally free and \(W\to H^0(E)\) is injective. Put \(r:=\rk E\).
Choosing \(m\) distinct points in the locally free locus gives
\(\dim W=mr\). For an arbitrary point \(x\), choose distinct points
\(x_2,\ldots,x_m\) in the locally free locus, all different from \(x\).
Evaluation on \(D=x+x_2+\cdots+x_m\) gives
\(mr=\dim W=\dim E_x+(m-1)r\), so \(\dim E_x=r\). On the smooth curve
\(\PP^1\), a torsion summand would increase the dimension of the fiber
at a point in its support. Thus \(E\) is torsion-free and hence locally
free. Finally, the composite \(W\to H^0(E)\to H^0(E|_D)\) is an
isomorphism, so its first arrow is injective.
\end{proof}

For a morphism \((\varphi,\psi):(E,W)\to(F,Z)\) and a reduced divisor
\(D\) of degree \(m\), the evaluation diagram
\[
\begin{tikzcd}[column sep=large,row sep=large]
W \arrow[r,"\psi"] \arrow[d,"\ev_D^E"',"\simeq"]
& Z \arrow[d,"\ev_D^F"',"\simeq"] \\
H^0(E|_D) \arrow[r,"H^0(\varphi|_D)"]
& H^0(F|_D)
\end{tikzcd}
\]
commutes. The vertical arrows are isomorphisms because both systems are
\(m\)-unisolvent.

\begin{lemma}\label{lem:constant-rank}
Let \((\varphi,\psi):(E,W)\to(F,Z)\) be a morphism in
\(\mathcal U_m\). Then \(\rk(\varphi_x)\) is independent of
\(x\in\PP^1\). Consequently, \(\ker\varphi\), \(\im\varphi\) and
\(\coker\varphi\) are vector bundles.
\end{lemma}

\begin{proof}
Let \(\rho\) be the generic rank of \(\varphi\). On a nonempty open
subset \(U\subset\PP^1\), one has \(\rk(\varphi_x)=\rho\), and
everywhere \(\rk(\varphi_x)\leq\rho\). Choosing \(m\) distinct points in
\(U\) and using the evaluation diagram gives \(\rk\psi=m\rho\). If
\(\rk(\varphi_x)<\rho\) at some point \(x\), choosing \(m-1\) further
points in \(U\) would instead give \(\rk\psi<m\rho\), a contradiction.
Thus \(\varphi\) has constant rank. Its kernel, image and cokernel are
therefore locally free.
\end{proof}

\begin{theorem}\label{thm:abelian}
The category \(\mathcal U_m\) is an extension-closed abelian subcategory
of \(\mathcal C\). Every object of \(\mathcal U_m\) has finite length.
\end{theorem}

\begin{proof}
The category \(\mathcal U_m\) is full and additive inside \(\mathcal C\).
By the standard criterion in
\cite[Remark~8.3.22 (ii)]{KashiwaraSchapira}, it is enough to show that
kernels and cokernels of its morphisms again lie in \(\mathcal U_m\).
Let \((\varphi,\psi):(E,W)\to(F,Z)\) be such a morphism and put
\(K:=\ker\varphi\) and \(C:=\coker\varphi\). These are vector bundles
by Lemma~\ref{lem:constant-rank}. For every reduced divisor \(D\) of
degree \(m\), restriction gives the exact sequence
\(0\to H^0(K|_D)\to H^0(E|_D)\to H^0(F|_D)\to H^0(C|_D)\to0\).
For the kernel, we compare the commutative ladder
\[
\begin{array}{ccccccc}
0 & \longrightarrow & \ker\psi & \longrightarrow & W & \longrightarrow & Z \\
  & & \big\downarrow & & \big\downarrow {\scriptstyle\sim}
  & & \big\downarrow {\scriptstyle\sim} \\
0 & \longrightarrow & H^0(K|_D) & \longrightarrow & H^0(E|_D)
  & \longrightarrow & H^0(F|_D).
\end{array}
\]
The two right vertical arrows are isomorphisms. Exactness therefore
identifies \(\ker\psi\) with \(H^0(K|_D)\). Similarly, the ladder
\[
\begin{array}{ccccccc}
W & \longrightarrow & Z & \longrightarrow & \coker\psi & \longrightarrow & 0 \\
\big\downarrow {\scriptstyle\sim} && \big\downarrow {\scriptstyle\sim}
&& \big\downarrow && \\
H^0(E|_D) & \longrightarrow & H^0(F|_D) & \longrightarrow & H^0(C|_D)
& \longrightarrow & 0
\end{array}
\]
identifies \(\coker\psi\) with \(H^0(C|_D)\). By
Lemma~\ref{lem:full-subcategory}, both \((K,\ker\psi)\) and
\((C,\coker\psi)\) belong to \(\mathcal U_m\). Thus \(\mathcal U_m\)
is abelian.

Now consider an exact sequence in \(\mathcal C\),
\(0\to(E_1,W_1)\to(E,W)\to(E_2,W_2)\to0\), whose two end terms belong
to \(\mathcal U_m\). The sheaf \(E\) is an extension of vector bundles
and is therefore locally free. For every reduced divisor \(D\) of
degree \(m\), there is a commutative diagram with exact rows
\[
\begin{array}{ccccccccc}
0&\longrightarrow&W_1&\longrightarrow&W&\longrightarrow&W_2&\longrightarrow&0\\
&&\big\downarrow {\scriptstyle\sim}&&\big\downarrow
&&\big\downarrow {\scriptstyle\sim}&&\\
0&\longrightarrow&H^0(E_1|_D)&\longrightarrow&H^0(E|_D)
&\longrightarrow&H^0(E_2|_D)&\longrightarrow&0.
\end{array}
\]
The two outer vertical arrows are isomorphisms, so the middle arrow is
also an isomorphism. Lemma~\ref{lem:full-subcategory} shows that
\((E,W)\in\mathcal U_m\). Hence \(\mathcal U_m\) is extension-closed.

Finally, the zero system is the only object of rank zero. A nonzero
proper subobject and its quotient both have positive rank, so rank
strictly increases along a strictly increasing chain of subobjects and
strictly decreases along a strictly decreasing chain. Either chain has
length at most the rank of the original object. Thus every object has
finite length.
\end{proof}

An object of \(\mathcal U_m\) is \emph{simple} if it is nonzero and has no
nonzero proper subobject. The theorem gives a finite Jordan--H\"older
filtration for every object. Its simple factors, counted with
multiplicity, are uniquely determined up to isomorphism and order.

\subsection{Simple Hermite systems and low-rank factors}

The Hermite systems of Example~\ref{par:split-hermite-models} are simple
in every rank. Together with the rank-two classification, this determines
all simple objects of ranks one and two.

\begin{proposition}\label{prop:hermite-simple}
The Hermite system \(H_{m,k}\) is a simple \(m\)-unisolvent system of
rank \(k\).
\end{proposition}

\begin{proof}
Put \(L:=\OO_{\PP^1}(mk-1)\), so
\(H_{m,k}=(J^{k-1}L,H^0(L))\), with sections included by their
\((k-1)\)-jets. A nonzero proper subobject
\((F,U)\subset H_{m,k}\) would have rank \(s\) with \(1\leq s<k\).
The jet-truncation morphism \(J^{k-1}L\to J^sL\) gives a composite
\(U\otimes\OO_{\PP^1}\to F\to J^{k-1}L\to J^sL\) whose image has rank
at most \(s\) at every point. Since \(\dim U=ms\geq s+1\), choose
linearly independent sections \(f_0,\ldots,f_s\in U\subset H^0(L)\).
On an affine coordinate chart with a local frame of \(L\), their
Wronskian \(\det(f_i^{(j)})_{0\leq i,j\leq s}\) vanishes identically.
In characteristic zero, the Wronskian criterion implies that the
\(f_i\) are linearly dependent over \(\CC\), a contradiction.
\end{proof}

\begin{corollary}\label{cor:low-rank-simple}
Up to isomorphism, \(L_m\) is the unique simple object of rank one and
\(H_{m,2}\) is the unique simple object of rank two in \(\mathcal U_m\).
\end{corollary}

\begin{proof}
Proposition~\ref{prop:low-rank} gives the unique rank-one object \(L_m\)
and the two rank-two objects \(L_m\oplus L_m\) and \(H_{m,2}\). A
rank-one object has no nonzero proper subobject. The split rank-two
object is not simple, whereas \(H_{m,2}\) is simple by
Proposition~\ref{prop:hermite-simple}.
\end{proof}

\subsection{Extensions and the diagonal exponent}

We record the extension sequence needed below. The extension-closedness
in Theorem~\ref{thm:abelian} identifies Yoneda extensions in
\(\mathcal U_m\) with extensions of coherent systems having the same
end terms.

\begin{proposition}\label{prop:ext-sequence}
Let \(S_1=(E_1,W_1)\) and \(S_2=(E_2,W_2)\) be objects of
\(\mathcal U_m\). There is a natural exact sequence
\[
\begin{aligned}
0&\to\Hom_{\mathcal U_m}(S_1,S_2)\to\Hom(E_1,E_2)
\to\Hom(W_1,H^0(E_2)/W_2)\\
&\to\Ext^1_{\mathcal U_m}(S_1,S_2)\to\Ext^1(E_1,E_2)
\to\Hom(W_1,H^1(E_2)).
\end{aligned}
\]
In particular, \(\Ext^1_{\mathcal U_m}(L_m,L_m)=0\).
\end{proposition}

\begin{proof}
By Theorem~\ref{thm:abelian}, the displayed sequence is the standard
extension sequence for coherent systems; see
\cite[Proposition~3.1]{LelliChiesa}. For \(S_1=S_2=L_m\), the chosen
subspace is the full space of global sections, and
\(\Ext^1(\OO_{\PP^1}(m-1),\OO_{\PP^1}(m-1))=H^1(\OO_{\PP^1})=0\).
The last assertion follows from exactness.
\end{proof}

The diagonal exponent also has a direct relation to extensions. For an
object \(S=(E,W)\), write \(e(S)\) for its diagonal exponent and put
\(a(S):=(e(S)-\rk E)/2\). By Theorem~\ref{thm:hilbert-discriminant},
\(a(S)\) is the multiplicity of the Hilbert evaluation determinant along
\(\Delta_m\). We set \(e(0)=a(0)=0\).

\begin{proposition}\label{prop:exponent-additivity}
For an exact sequence \(0\to S_1\to S\to S_2\to0\) in
\(\mathcal U_m\), one has
\(e(S)=e(S_1)+e(S_2)\) and \(a(S)=a(S_1)+a(S_2)\).
If \(Q_1,\ldots,Q_\ell\) are the Jordan--H\"older factors of \(S\), then
\(e(S)=\sum_{i=1}^{\ell}e(Q_i)\) and
\(a(S)=\sum_{i=1}^{\ell}a(Q_i)\).
\end{proposition}

\begin{proof}
Write \(S_i=(E_i,W_i)\) and \(S=(E,W)\). By
Proposition~\ref{prop:detcriterion},
\(e(S)(m-1)=\deg E=\deg E_1+\deg E_2
=(e(S_1)+e(S_2))(m-1)\). Since \(m\geq2\), this proves the first
identity. Ranks are also additive, so the identity for \(a\) follows.
Applying the two identities successively to a Jordan--H\"older
filtration gives the last assertions.
\end{proof}

\subsection{The rank-three threshold}

The possible low-rank factors and the additivity of the exponent give a
numerical restriction on every nonsimple rank-three system.

\begin{proposition}\label{prop:rank-three-simple}
The diagonal exponent of a nonsimple rank-three \(m\)-unisolvent system
belongs to \(\{3,5\}\). Every rank-three system of exponent \(3\) is
isomorphic to \(L_m^{\oplus3}\). Consequently, every rank-three
\(m\)-unisolvent system of exponent greater than \(5\) is simple.
\end{proposition}

\begin{proof}
Choose a nonzero proper subobject \((F,Z)\subset(E,W)\) and denote its
quotient by \((Q,Y)\). The two ranks are \(1\) and \(2\) in some order.
By Proposition~\ref{prop:low-rank}, the rank-one term has exponent \(1\)
and the rank-two term has exponent \(2\) or \(4\). Hence
Proposition~\ref{prop:exponent-additivity} gives
\(e(E,W)=e(F,Z)+e(Q,Y)\in\{3,5\}\). Independently of simplicity,
Proposition~\ref{prop:minimal-exponent} identifies every rank-three
system of exponent \(3\) with \(L_m^{\oplus3}\). The final assertion is
the contrapositive of the first.
\end{proof}

The remaining boundary value is rigid when \(m=2\).

\begin{proposition}\label{prop:rank-three-exponent-five}
Every rank-three \(2\)-unisolvent system of diagonal exponent \(5\) is
isomorphic to \(H_{2,2}\oplus L_2\). In particular, no such system is
simple.
\end{proposition}

\begin{proof}
Let \((E,W)\) be such a system. Proposition~\ref{prop:detcriterion} gives
\(\deg E=5\), while Lemma~\ref{lem:unisolvent-splitting-lower-bound}
shows that every summand in the Birkhoff--Grothendieck decomposition of
\(E\) has degree at least one. After ordering the summands, the only
possibilities are \(E\simeq\OO_{\PP^1}(3)\oplus\OO_{\PP^1}(1)^{\oplus2}\)
or \(E\simeq\OO_{\PP^1}(2)^{\oplus2}\oplus\OO_{\PP^1}(1)\).

Suppose first that the splitting type is \((3,1,1)\). Let
\(F\simeq\OO_{\PP^1}(3)\) be the degree-three summand and put
\(Z:=W\cap H^0(F)\). Since \(h^0(E)=8\), \(h^0(F)=4\) and
\(\dim W=6\), one has \(\dim Z\geq2\). For every reduced divisor \(D\)
of degree two, the map \(Z\to H^0(F|_D)\) is injective: its composite
with \(H^0(F|_D)\to H^0(E|_D)\) is the restriction of the evaluation
isomorphism for \((E,W)\). Its target has dimension two, so
\(\dim Z=2\) and \((F,Z)\) is a rank-one \(2\)-unisolvent system. This
contradicts Proposition~\ref{prop:low-rank}, since its line bundle is
\(\OO_{\PP^1}(3)\), not \(\OO_{\PP^1}(1)\).

Thus \(E\simeq\OO_{\PP^1}(2)^{\oplus2}\oplus\OO_{\PP^1}(1)\). Put
\(F:=\OO_{\PP^1}(2)^{\oplus2}\) and \(Z:=W\cap H^0(F)\). The same
argument gives \(\dim Z\geq6+6-8=4\), while evaluation embeds \(Z\)
into the four-dimensional space \(H^0(F|_D)\) for every reduced divisor
\(D\) of degree two. Hence \(\dim Z=4\), and \((F,Z)\) is a rank-two
\(2\)-unisolvent system of degree four. Proposition~\ref{prop:low-rank}
identifies it with \(H_{2,2}\). Its quotient in \(\mathcal U_2\) has
rank one and is therefore \(L_2\). We obtain an exact sequence
\(0\to H_{2,2}\to(E,W)\to L_2\to0\).

We claim that \(\Ext^1_{\mathcal U_2}(L_2,H_{2,2})=0\). Write
\(H_{2,2}=(\OO_{\PP^1}(2)^{\oplus2},W_H)\). The two objects \(L_2\)
and \(H_{2,2}\) are simple and have different ranks, so
\(\Hom_{\mathcal U_2}(L_2,H_{2,2})=0\). Proposition~\ref{prop:ext-sequence}
therefore gives an injection
\[
\Hom\bigl(\OO_{\PP^1}(1),\OO_{\PP^1}(2)^{\oplus2}\bigr)
\longrightarrow
\Hom\bigl(H^0(\OO_{\PP^1}(1)),H^0(\OO_{\PP^1}(2)^{\oplus2})/W_H\bigr).
\]
Both spaces have dimension four, so this map is an isomorphism. Moreover,
one has
\[
\Ext^1(\OO_{\PP^1}(1),\OO_{\PP^1}(2)^{\oplus2})
=H^1(\OO_{\PP^1}(1))^{\oplus2}=0.
\]
Exactness gives the claim.
Consequently, the preceding short exact sequence splits, and
\((E,W)\simeq H_{2,2}\oplus L_2\).
\end{proof}

\iffalse For \(m=2\), these two boundary values also have concrete Grassmannian descriptions.
For exponent \(3\), Proposition~\ref{prop:minimal-exponent} gives the sum
of three moving lines, each parametrizing \(\PP^1\) in a fixed
two-dimensional summand of \(V\). For exponent \(5\), there is a
decomposition \(V=A\oplus B\), with \(\dim A=4\) and \(\dim B=2\), such
that \(U_t=\widehat T_{c(t)}C\oplus\ell(t)\). Here
\(c:\PP^1\to\PP(A)\) is a twisted cubic with image \(C\), and
\(t\mapsto\ell(t)\in\PP(B)\) is a degree-one rational normal embedding.
Thus the exponent-five model combines the tangent-line family of a
twisted cubic with the rank-one Vandermonde curve.\fi

The next example gives a simple system at the next possible exponent.

\begin{example}[A simple exponent-seven system]\label{ex:rank-three-exponent-seven}
Let \(V=\CC^6\), with standard basis \(e_1,\ldots,e_6\). On the affine
chart \(\Aff^1\subset\PP^1\), consider the polynomial row frame
\[
R_X(t)=
\begin{pmatrix}
1&t&t^2&t^3&t^4&0\\
0&1&2t&3t^2&4t^3&0\\
0&0&1&3t&6t^2&t
\end{pmatrix}.
\]
Let \(U_t\) be its row space. A direct calculation gives
\(\det\begin{pmatrix}R_X(s)\\R_X(t)\end{pmatrix}=-(s-t)^7\), so
\(U_s\oplus U_t=V\) whenever \(s\neq t\) are finite.

We check the point at infinity. The three-by-three minors of \(R_X(t)\)
have degree at most \(7\). The minors of the first and last three columns
are \(p_{123}(t)=1\) and \(p_{456}(t)=t^7\), and the latter is the only
minor of degree \(7\). Their degree-seven homogenizations therefore have
no common zero and extend the Pl\"ucker map to \(\PP^1\). Since the
Grassmannian is closed in its Pl\"ucker embedding, the extension still
takes values in \(\Gr(3,V)\), with
\(U_\infty=\Span\{e_4,e_5,e_6\}\). Moreover,
\[
\det\begin{pmatrix}R_X(t)\\ 0_{3\times3}\quad I_3\end{pmatrix}
=\det\begin{pmatrix}1&t&t^2\\0&1&2t\\0&0&1\end{pmatrix}=1.
\]
Thus \(U_t\oplus U_\infty=V\) for every finite \(t\). The resulting
morphism \(\PP^1\to\Gr(3,6)\) is \(2\)-direct-sum. Its associated
system has diagonal exponent \(7\), by the two-point determinant above,
and is simple by Proposition~\ref{prop:rank-three-simple}.
\end{example}

\begin{corollary}\label{cor:rank-three-threshold}
The least diagonal exponent of a simple rank-three \(2\)-unisolvent
system is \(7\).
\end{corollary}

\begin{proof}
Corollary~\ref{cor:exponent-bounds} gives \(e\geq3\) and
\(e\equiv1\pmod2\). The exponent-three system is \(L_2^{\oplus3}\),
and every exponent-five system is \(H_{2,2}\oplus L_2\), by
Propositions~\ref{prop:rank-three-simple}
and~\ref{prop:rank-three-exponent-five}. Neither is simple.
Example~\ref{ex:rank-three-exponent-seven} gives a simple system of
exponent \(7\), proving the assertion.
\end{proof}

\section*{Acknowledgements}
\phantomsection
\addcontentsline{toc}{section}{Acknowledgements}
The authors thank Professor Baohua Fu and the Institute of Mathematics of the 
Chinese Academy of Sciences for organizing the "Workshop on Symplectic Varieties and Projective Geometry" and for their
warm hospitality in Beijing, which made possible the discussions that led
to this collaboration.

\bigskip
\noindent\begin{minipage}{\textwidth}
\small
\textsc{Duo Li}\\
School of Mathematics, Sun-Yat Sen University,\\
Zhuhai, Guangdong 519082, PR China\\
\textit{Email:} \href{mailto:liduo5@mail.sysu.edu.cn}{\texttt{liduo5@mail.sysu.edu.cn}}
\end{minipage}

\medskip
\noindent\begin{minipage}{\textwidth}
\small
\textsc{Yanjie Li}\\
Academy of Mathematics and Systems Science, Chinese Academy of Sciences,\\
55 Zhong Guan Cun East Road, Beijing 100190, China;\\
University of Chinese Academy of Sciences, Beijing, China\\
\textit{Email:} \href{mailto:liyanjie241@mails.ucas.ac.cn}{\texttt{liyanjie241@mails.ucas.ac.cn}}
\end{minipage}

\medskip
\noindent\begin{minipage}{\textwidth}
\small
\textsc{Francesco Russo}\\
Dipartimento di Matematica e Informatica, Università di Catania,  \\ Città Universitaria, Viale A. Doria 6,
95125 Catania, Italy\\
\textit{Email:} \href{mailto:francesco.russo@unict.it}{\texttt{francesco.russo@unict.it}}
\end{minipage}

\end{document}